\documentclass[oneside,reqno]{amsart}
\usepackage{amsmath}
\usepackage{amssymb,latexsym,amstext,amsthm,verbatim}
\usepackage{eucal} 
\usepackage{mathdots} 

\theoremstyle{plain}
\newtheorem{theorem}{Theorem}[section]
\newtheorem{lemma}[theorem]{Lemma}
\newtheorem{corollary}[theorem]{Corollary}
\newtheorem{proposition}[theorem]{Proposition}

\theoremstyle{definition}
\newtheorem{definition}[theorem]{Definition}

\newtheorem{example}[theorem]{Example}
\newtheorem{remark}[theorem]{Remark}

\newcommand\os{\overset}

\newcommand\bX{{\mathbb X_k}}
\newcommand\bV{{\mathbb V_k}}
\newcommand\cA{\mathcal A}

\newcommand\cD{{\mathcal D}}
\newcommand\cI{{\mathcal I}}
\newcommand\cJ{{\mathcal J}}
\newcommand\cL{{\mathcal L}}

\newcommand\cX{{\mathcal X}}

\newcommand\bi{{\mathbf i}}
\newcommand\bj{{\mathbf j}}
\newcommand\bq{{\mathbf q}}
\newcommand\bbZ{\mathbb Z}

\newcommand\andd{\quad \text{and} \quad}
\newcommand\bC{\bar C}
\newcommand\centA{C(\cA)}
\newcommand\centL{C(\cL)}
\newcommand\modi{{\bar{\imath}}}
\newcommand\modj{{\bar{\jmath}}}

\newcommand\Mat{M}
\newcommand\Mloop{M}
\newcommand\ot{\otimes}
\newcommand\otk{\otimes_k}
\newcommand\q{\quad}
\newcommand \set[1]{\{ \, #1 \, \}}
\newcommand \sll{\operatorname{sl}_{\ell+1}(k)}
\newcommand\Stn{S^{\otimes n}}
\newcommand\Stwo{S^{\otimes 2}}
\newcommand\Stp{S^{\otimes p}}
\newcommand\Stpplus{S^{\otimes (p+1)}}
\newcommand\strange{k[u_1,u_2^{\pm 1},w]_\rho}
\newcommand\strangep{k[u_1,u_2^{\pm 1},w]_{\rho'}}
\newcommand\suchthat{\mid }
\newcommand\tA{\wt\cA}
\newcommand\tcentA{\wt{C(\cA)}}
\newcommand\thsup{\text{th}}
\newcommand\us{\underset}

\newcommand\wt{\widetilde}
\newcommand\Zmod[1]{\bbZ_{#1}}

\newcommand\Aut{\operatorname{Aut}}
\newcommand\charr{\operatorname{char}}
\newcommand\diag{\operatorname{diag}}
\newcommand\Dim{\operatorname{Dim}} 
\newcommand\GL{\operatorname{GL}}

\newcommand\End{\operatorname{End}}
\newcommand\Hom{\operatorname{Hom}}
\newcommand\Mult{\operatorname{Mult}}
\newcommand\one{\operatorname{1}}
\newcommand\spann{\operatorname{span}}
\newcommand \speciallinear{\operatorname{sl}}

\begin{document}
\setcounter{page}{1}

\title{Iterated Loop Algebras}
\author{Bruce Allison}
\address[Bruce Allison]
{Department of Mathematical and Statistical Sciences\\ University of
Alberta\\Edmonton, Alberta, Canada T6G 2G1}
\email{ballison@math.ualberta.ca}
\author{Stephen Berman}
\address[Stephen Berman]
{Department of Mathematics and Statistics\\
University of Sask\-at\-che\-wan\\ Saskatoon, Saskatchewan, Canada S7N 5E6}
\email{berman@math.usask.ca}
\author{Arturo Pianzola}
\address[Arturo Pianzola]
{Department of Mathematical and Statistical Sciences\\ University of
Alberta\\Edmonton, Alberta, Canada T6G 2G1}
\email{a.pianzola@ualberta.ca}
\thanks{The authors gratefully acknowledge the support
of the Natural Sciences and Engineering Research Council of Canada.}
\date{\today}

\begin{abstract}
Iterated loop algebras are by
definition obtained by repeatedly applying
the loop construction, familiar from the theory of affine Kac-Moody
Lie algebras, to a given base algebra.  Our interest
in this iterated construction is motivated by its use in the realization
of extended affine Lie algebras, but the construction also appears
naturally in the study of other classes of algebras.
This paper consists of a detailed study of the basic properties
of iterated loop algebras.
\end{abstract}

\maketitle

\section{Introduction}
Over  the past 35 years affine  Kac-Moody Lie algebras have been at the centre of a considerable amount of beautiful
mathematics and theoretical physics. As of late, and perhaps
influenced by some of the newest theories in physics, the need
seems to have arisen for some ``higher nullity'' generalizations
of affine Kac-Moody Lie algebras.  It is still too early to decide
what the correct final choice for these algebras will be, but it
is fair to say notwithstanding, that Lie algebras graded by root
systems and extended affine Lie algebras (EALAs) will play a
prominent role in the  process \cite{BM,BZ,AABGP,SY}.

Recall that given a
$\Zmod{m}$-grading $\Sigma = \{\cA_\modi\}_{\modi \in \Zmod{m}}$
of an algebra $\cA$ over a field $k$, the \emph{loop algebra} of $\Sigma$
based on $\cA$
is the subalgebra
\[L(\cA,\Sigma) := \bigoplus_{i\in \bbZ} \cA_\modi \otk z^i\]
of $\cA\otk k[z,z^{-1}]$.
Using this beautiful construction, V.~Kac showed that
(the derived algebra modulo its centre of) any complex
affine Kac-Moody
Lie algebras can be obtained as a loop algebra of a finite
dimensional simple Lie algebra \cite{K1}.
The loop construction makes it
clear, among other things, that the affine algebras are objects of
nullity one in a sense that can be made precise.  Indeed, in
EALA theory, where the concept of nullity is well-defined, one
finds that finite dimensional simple algebras are precisely the
(tame) EALAs of nullity zero whereas affine algebras are
precisely the (tame) EALAs of
nullity one \cite{ABGP}.

It thus seems almost inevitable to ask whether, starting from an affine Kac-Moody
Lie algebra and applying the loop construction, one obtains an extended
affine Lie algebra of nullity 2.  This and related questions have been
investigated in some detail in \cite{W,Po,ABP1,ABP2,vdL}.
In our work on this topic, as well as in \cite{vdL},  it became clear that some
advantages are to be had by thinking of loop
algebras based on an affine algebra as being obtained from a finite
dimensional Lie algebra by applying the loop construction
\textit{twice} (the advantages stemming from the fact that in this
case the ``bottom'' algebra, namely the finite dimensional one,
is much simpler than the affine algebra).  As the  reader
will have surmised by now, the study of these ``iterated loop
algebras'' took on a life of its own and became the subject of the
present paper.

In general, if $\cA$ is an (arbitrary) algebra over $k$, an
\emph{$n$-step iterated loop algebra} based on $\cA$
is an algebra that can be obtained starting from $\cA$ by a sequence of
$n$ loop constructions, each based on the algebra obtained at the previous
step (see Definition \ref{def:iteratedloop}).  Far from being a mere generalization of the loop
construction, iterated loop algebras seem to yield interesting
mathematical objects in a natural way.  Even when the resulting
objects are known, the new point of view can be illuminating.
As an example, we see in  Example~\ref{ex:quantum} that algebras
representing elements of the Brauer group of the
ring $k[t_1^{\pm 1},t_2^{\pm 2}]$ are
obtained as 2-step iterated loop algebras of $M_n(k)$.
This information is not apparent if one thinks in terms of
single loop algebras of $M_\ell(k[t_1^{\pm 1}])$.

This paper contains a detailed study of the basic properties
of iterated loop algebras.  We begin in Section \ref{sec:gen}
by recording some simple properties of the centroid of an algebra.
In the rest of Section \ref{sec:gen} and in Section \ref{sec:primepfgc}
we define and give the basic properties
of a very important class of algebras which for lack of a better name
we have simply referred to as
pfgc algebras (nonzero, perfect, and finitely generated as modules over their
centroids).  The property of being a prime pfgc algebra arises naturally in the study
of iterated loop algebras since this property is carried over to a loop algebra
(and hence to an iterated loop algebra) from its base.  In contrast
the property of finite dimensional central simplicity
certainly does not carry over in the same way.  After this discussion
of pfgc algebras we establish in Section \ref{sec:loop} some basic properties of
(one step) loop algebras.

The main results of the paper appear in Sections \ref{sec:iterated},
\ref{sec:centroiditerated}, \ref{sec:untwist} and \ref{sec:permanence} .
These all deal with properties of
an $n$-step iterated loop algebra $\cL$ based on a pfgc algebra $\cA$.
First Theorem \ref{thm:basic} establishes
a long list of properties that carry over from $\cA$ to $\cL$.
In particular, it is shown (as mentioned above)
that if $\cA$ is a prime pfgc algebra then so is
$\cL$.
Next Theorem \ref{thm:centloopn} shows that the centroid $C(\cL)$ of
$\cL$ is itself
an $n$-step iterated loop algebra of the centroid of $\cA$.  The same theorem describes
a method of calculating $C(\cL)$ explicitly.  Then
Theorem \ref{thm:untwist} shows that $\cL$ can be ``untwisted''
by a base ring extension of  $C(\cL)$ that is free of finite rank.
That
is, the algebra $\cL$ (after such a base ring extension) becomes isomorphic
to the iterated loop algebra obtained using only the trivial gradings at each stage.
Finally, Section \ref{sec:permanence} deals with the concept of
type of an algebra (which is motivated by
the concept of   type in terms of root
systems which exists in Lie theory).  The main result, Theorem
\ref{thm:permanence}, states that   type
cannot change under the loop construction.

Each of the main results in Sections
\ref{sec:centroiditerated}, \ref{sec:untwist} and
\ref{sec:permanence} has several corollaries that are discussed in
the respective sections.  To give  one important example, we show
in Section \ref{sec:permanence} that if $\cL$ is an $n$-step iterated loop
algebra based on a finite dimensional split simple Lie algebra $\cA$
over a field of characteristic 0 then both $\cA$
and $n$ are isomorphism invariants of $\cL$ (see Corollary \ref{cor:Liebase}).
This  result will play a crucial role in
our forthcoming work on the classification of the centreless cores
of EALA's of nullity 2 \cite{ABP3}.

In the last section, Section \ref{sec:twostep},
we look closely at 2-step iterated loop algebras.  If the base  algebra
is finite dimensional and central simple, these 2-step iterated loop algebras
come in two kinds,
depending on the structure of their centroids.  We illustrate
this fact along with many of the concepts discussed in the paper by describing
two examples dealing respectively
with Lie algebras and associative   algebras.

\section{Centroids and pfgc algebras}
\label{sec:gen}

In this section we record some basic facts about centroids, and we
define a class of algebras, which we call pfgc algebras, that will
play an important role in the study of loop algebras.  A good
basic reference on the centroid is   \cite[Ch.~X, \S~1]{J1}.

\bigskip

{\bf Terminology and notation:}
A  \textit{ring} will mean a unital commutative associative ring.
Homomorphisms, subrings and modules for rings will always be assumed to be unital.

A  \textit{base change} will mean a homomorphism $\upsilon: B \to B'$ of rings.
This base change is said to be \textit{free}
(respectively \textit{flat, faithfully flat}) if $B'$ is a free (respectively flat,
faithfully flat) $B$-module.  Note that if $\upsilon: B \to B'$ is free and $B'\ne 0$,
then $\upsilon$ is faithfully flat and hence flat \cite[\S~I.3.1, Example~(2)]{B:CA}.
An injective base change $\upsilon : B \to B'$ will be called an
\textit{extension of rings}, in which case we often identify
$B$ as a subring of $B'$ and denote the extension by  $B'/B$.

If $B$ is a ring, an \textit{algebra}
over $B$ will mean  a $B$-module  $\cA$ together with a $B$-bilinear product (which is not necessarily
associative, commutative or unital).  If $\cA$ and $\cA'$ are $B$-algebras,
we use the notation
\[\cA \simeq_B \cA'\]
to mean that $\cA$ and $\cA'$ are isomorphic as $B$-algebras.
If $\cA$ is an algebra over $B$ and
$\upsilon : B\to B'$ is a base change,
we will denote by $\cA \ot_B B'$ the (unique) $B'$-algebra
which is obtained from $\cA$ by base change \cite[Ch.~III,
\S~1.5]{B:Alg}.

\bigskip

\textit{For the rest of the section we assume that $B$ is
a ring, and that $\cA$ is
a $B$-algebra.}   Note that $\cA$ can also be regarded as
$\bbZ$-algebra under the natural action of $\bbZ$ on $\cA$.

We  now recall the definition of the centroid of $\cA$
\cite[Ch.~X, \S~1]{J1}.

\begin{definition}
\label{def:centroid} (i) For $a\in{\cA}$ consider the two maps
from ${\cA}$ to ${\cA}$
$$
a_L : x\mapsto ax \q{\rm and}\q a_R:x\mapsto xa.
$$
The  \textit{multiplication algebra}  of $\cA$ \cite[Ch.~X,
\S~1]{J1} is defined to be the $B$-subalgebra $\Mult_B(\cA)$ of
$\End_B(\cA)$ generated by $\set{\one}\cup\set{a_L \suchthat a\in
\cA}\cup\set{a_R \suchthat a\in \cA}$.

(ii)  The set $C_B(\cA)$ of elements of $\End_B(\cA)$ that commute
with the action of $\Mult_B(\cA)$  is called the \emph{centroid}
of $\cA$.
Equivalently
\[
C_B(\cA) := \{\chi  \in \;\End_B(\cA): \chi(xy) = \chi(x)y
=x\chi(y)\q \text{\rm for all} \q x,y\in \cA\}.
\]
(The notation $\operatorname{Cent}_B(\cA)$ has been used for the centroid
in some articles, for example in \cite{ABP2}.  We are using the abbreviated
notation $C_B(\cA)$ since it will arise frequently.)
Clearly $C_B(\cA)$ is a $B$-subalgebra of $\End_B(\cA),$ and therefore
$\cA $ can be viewed in a natural way as a left $C_B(\cA)$-module
by defining $\chi  \cdot x=\chi(x).$

(iii) For $b\in B$ we define $\lambda_\cA(b)\in \End_B(\cA)$ by
\[
\big(\lambda_{\cA}(b)\big)(x) = b\cdot x.
\]
Clearly $\lambda_{\cA}(b)\in  C_B(\cA)$ since $\cA$ is a
$B$-algebra.  Then the map $\lambda_{\cA} : B \to C_B(\cA)$ is a
ring homomorphism, and $C_B(\cA)$ \textit{is a unital associative
$B$-algebra via this map}. Furthermore, if $\cA$ is  a faithful
$B$-module then $B$ can be identified with a subring of (the
centre of) the centroid $C_B(\cA).$

(iv) The $B$-algebra $\cA$ is said to be  \emph{central}
(or central over $B$) if
$\lambda_{\cA} : B \to C_B(\cA)$ is an isomorphism.

(v) The  \textit{centre}  of $\cA$ is defined to be the
set $Z(\cA)$ of elements in $\cA$ that commute and
associate with all elements of $\cA$.
Then $Z(\cA)$ is a $B$-subalgebra of $\cA$.
If $\cA$ is unital,
the map which sends $z$ to left multiplication by $z$
is a $B$-algebra isomorphism of $Z(\cA)$ onto $C_B(\cA)$
\cite[\S~1]{EMO}.
\end{definition}

The following is clear:

\begin{lemma}
\label{lem:centinduce} Suppose that $\cA$ and $\cA'$ are
$B$-algebras and  $\rho  :\cA\to \cA'$ is a $B$-algebra
isomorphism.  Then $\rho$ induces a $B$-algebra isomorphism
$C_B(\rho) : C_B(\cA)\to C_B(\cA')$ defined by $\chi  \mapsto \rho
\chi  \rho^{-1}$.
\end{lemma}

The formation of the centroid does not commute with base change.
Nonetheless these two processes do commute in two important cases
that we now describe. If $B\to B'$ is a homomorphism of rings, we define
\[\nu = \nu_{\cA,B,B'} : C_B(\cA)\ot_B B' \to
C_{B'}(\cA\ot_B B')\] to be the restriction of the canonical  map
$\End_B(\cA) \ot_B B' \to \;\End_{B'}(\cA \ot_B B')$. Then $\nu$
is a homomorphism, said to be canonical, of unital associative
$B'$-algebras.

\begin{lemma}
\label{lem:BaseChange} Suppose that $B\to B'$ is a homomorphism of
rings.  Then the map
$\nu_{\cA,B,B'} :C_B(\cA)\ot_B B' \to C_{B'}(\cA\ot_B B') $ is an
isomorphism of $B'$-algebras in the following  cases:
\begin{itemize}
\item[(a)] $\cA$ is finitely generated as a module over
its multiplication algebra $\Mult_B(\cA)$ and $B'$ is a free
$B$-module
\item[(b)] $B'$ is a finitely
generated projective $B$-module
\end{itemize}
\end{lemma}

\begin{proof}
(a)   Let $\{s_i\}_{i\in
I}$ be a basis of the $B$-module $B'$.

It is clear that $\nu$ is injective. Indeed  if $\sum \chi_i\ot s_i$ is
in the kernel of $\nu$ then $\sum \chi_i(x)\ot s_i =0$ for all $x$
in $\cA$ and so $\chi_i = 0$ for all $i$ in $I.$

To see that $\nu$ is onto, let $\chi  \in C_{B'}(\cA \ot_B B')$.
Then for
$x\in \cA$ we can write
$\chi(x\ot_B \one_{B'})$ uniquely as
\[
\chi(x\ot \one_{B'}) =\sum \chi_i(x)\ot s_i,
\]
where $\chi_i(x)\in \cA$ and only finitely many of these are
nonzero.  It is easy to see that for all $i\in I$ the map $\chi
_i:\cA\to\cA$ given by $\chi_i:x\mapsto \chi_i(x)$ is an element
of $C_B(\cA)$.  Thus to see that $\chi$ is an image under $\nu$ it
suffices to show that only finitely many of the maps $\chi_i$ are
nonzero.  For this let $\set{x_1,\dots,x_n}$ be a set of
generators of $\cA$ as a $\Mult_B(\cA)$-module. Then whenever
$\chi_i$ vanishes on all $x_j$'s we have
\[
\chi_i(\cA) =\chi_i\Big(\sum^n_{j=1} \Mult_B(\cA)\cdot x_j\Big) =
\sum^n_{j=1} \Mult_B(\cA)\cdot \chi_i(x_j)=0.
\]

(b) Consider the unique $B$-module homomorphism
\[\varphi_{B,\cA} : \End_B(\cA) \to \Hom_B(\cA\ot_B \cA,
\cA\oplus \cA)\] satisfying
\[
\varphi_{B,\cA}(f)(a_1\ot_B a_2) =\big(f(a_1a_2)
-f(a_1)a_2,f(a_1a_2)-a_1f(a_2)\big).
\]
By definition
\[
\ker(\varphi_{B,\cA}) =C_B(\cA).
\]
Also,  by standard properties of projective modules we obtain the
diagram
\smallskip
\[
\begin{smallmatrix}
0 \to &C_B(\cA)\ot_B B' &\to &\End_B(\cA)\ot_B B'
&\to &\Hom_B(\cA\ot_B\cA,\cA\oplus \cA) \ot_B B'\\
&\downarrow\nu   &&\Vert &&\Vert\\
0 \to &C_{B'}(\cA\ot_B B') &\to &\End_{B'}(\cA\ot_B B') &\to
&\Hom_{B'}\big((\cA\ot_B B')\ot_{B'} (\cA\ot_B B'), &\cA\ot_B B'
\oplus \cA\ot_B B'\big),
\end{smallmatrix}
\]
\smallskip
where the horizontal rows are exact.  Indeed the exactness of the
top row is by flatness of  the $B$-module $B'$ (every projective is flat).  The two
vertical isomorphisms come from $B'$ being a finitely generated
$B$-module which is projective \cite[Ch.~II, \S~5.3,
Prop.~7]{B:Alg}. It follows that $\nu$ is an isomorphism.
\end{proof}

The following important fact is proved  in \cite[Ch.~X, \S~1,
Theorem~3]{J1}:

\begin{lemma}
\label{lem:centalsimple} Suppose that $B$ is a field and $\cA$ is
finite dimensional and central simple over $B$. If $B'/B$ is a
field extension, then $\cA\ot_B B'$ is finite dimensional and
central simple over $B'$.
\end{lemma}

Next we consider gradings on $C_B(A)$ that are induced by gradings
on $\cA$. For this suppose that $\cA$ is $Q$-graded algebra over
$B$ where $Q$ is a finite abelian group.  Thus
\[
\cA = \bigoplus_{\alpha  \in Q} \cA_\alpha
\]
for some $B$-submodules $\cA_\alpha$ and $\cA_\alpha \cA_\beta
\subset \cA_{\alpha+\beta}.$ Then, since $Q$ is finite,
\[\End_B(\cA) = \bigoplus_{\lambda  \in Q} \End_B(\cA)_\lambda\]
is also a $Q$-graded $B$-algebra, where
\[
\End_B(\cA)_\lambda  =\{\theta  \in \End_B(\cA) \suchthat
\theta(\cA_\alpha  )\subset \cA_{\lambda  +\alpha} \text{ for all
} \alpha  \in Q\}.
\]
It is easy to check that $C_B(\cA)$ is a $Q$-graded $B$-subalgebra
of $\End_B(\cA)$, and so we have:

\begin{lemma}
\label{lem:graded} Suppose that $\cA$ is $Q$-graded algebra over
$B$ where $Q$ is a finite group. Then
\[
C_B(\cA) = \bigoplus_{\alpha  \in Q} C_B(\cA)_\lambda
\]
is a $Q$-graded algebra over $B$, where $C_B(\cA)_\lambda =
C_B(\cA) \cap \End_B(\cA)_\lambda$ for  all $\lambda \in Q.$
\end{lemma}

\begin{definition}
If ${\cI}$ and ${\cJ}$ are ideals of the $B$-algebra $\cA$ we
define
$$
{\cI}{\cJ} =\big\{\sum   x_iy_i : x_i \in {\cI},\,
y_i\in{\cJ}\big\}
$$
(finite sums of course).  Note that in general ${\cI}{\cJ}$ is not
an ideal of $\cA.$ We say that $\cA$ is \emph{perfect} if $\cA\cA
= \cA$.
\end{definition}

\begin{remark}
\label{rem:perfectinvariant} It is clear that $\cA$ is perfect as
a $B$-algebra if and only if $\cA$ is perfect as a $\bbZ$-algebra.
\end{remark}

\begin{lemma}
\label{lem:centroidinvariant} Assume $\cA$ is perfect. Then
\begin{itemize}
\item[(i)] $C_B(\cA)$ is commutative.
\item[(ii)]  $C_B(\cA) =
C_{\bbZ}(\cA).$
\end{itemize}
\end{lemma}

\begin{proof} (i) See  \cite[Ch.~X, \S~1, Lemma~1]{J1}.

(ii) We must show that any element $\chi\in C_\bbZ(\cA)$ is
$B$-linear. Indeed if $x,y\in \cA$ and $b\in B$ we have
$\chi(b\cdot(xy)) =\chi(x(b\cdot y)) = \chi(x)(b\cdot y) = b\cdot
(\chi(x)y) = b\cdot\chi(xy)$.
\end{proof}

We now  introduce a convenient acronym, pfgc, that will be used
throughout the paper.

\begin{definition}
\label{def:pfgc} A $B$-algebra $\cA$ is said to be \emph{pfgc} if
it satisfies the following conditions
\begin{itemize}
\item[P0.] $\cA\ne (0)$
\item[P1.] $\cA$ is perfect
\item[P2.] $\cA$ is finitely generated as a module over its
centroid $C_B(\cA).$
\end{itemize}
\end{definition}

\begin{remark}
\label{rem:pfgcinvariant} The notion of pfgc algebra $\cA$ is
independent of the base ring under which $\cA$ is viewed as an
algebra. More precisely, if $\cA$ is an algebra over $B$, it
follows from Remark \ref{rem:perfectinvariant} and
Lemma~\ref{lem:centroidinvariant}(ii) that $\cA$ is a pfgc algebra
over $B$ if and only if $\cA$ is a pfgc algebra over $\bbZ$
\end{remark}

We now summarize the basic facts that we will need about pfgc
algebras.

\begin{proposition}
\label{prop:fg}  Suppose that $\cA$ is a pfgc algebra over $B$.
Then
\begin{itemize}
    \item[(i)] $C_B(\cA)$ is
    a nonzero unital commutative associative $B$-algebra.
    \item[(ii)]  $\cA$ is finitely generated as a module over its
    multiplication algebra~$\Mult_B(\cA)$.
\end{itemize}
\end{proposition}

\begin{proof}
(i) Since $\cA$ is perfect and nonzero, this follows from
Lemma~\ref{lem:centroidinvariant}(i).

(ii) Let $C=C_B(\cA).$  Let $\{x_1,\dots,x_n\}\in \cA$ be such
that $\cA=\sum Cx_i$. For each $i$ we can write $x_i=\sum_j
y_{ij}z_{ij}$ (finite sum) for some $y_{ij}$ and $z_{ij}$ in
$\cA.$   Then
$$
\cA = \sum_{i} Cx_i =\sum_{i,j} C(y_{ij}z_{ij})
=\sum_{i,j}(Cy_{ij})z_{ij} \subset \sum_{i,j} \Mult_B(\cA)\cdot
z_{ij}  ,
$$
which shows that  $\cA$ is generated by the $z_{ij}$'s as an
$\Mult_B(\cA)$-module.
\end{proof}

\section{Prime pfgc algebras}
\label{sec:primepfgc}

In this section, we recall some basic facts about prime algebras
and consider in particular properties of prime pfgc algebras. A
good basic reference on prime nonassociative algebras and their
centroids is \cite{EMO}.

We suppose again in this section that $B$ is a ring and that $\cA$ is $B$-algebra.

\begin{definition}
The $B$-algebra $\cA$ is said to be \emph{prime} if for all ideals
${\cI}$ and ${\cJ}$ of the $B$-algebra $\cA$ we have
$$
{\cI}{\cJ} =0\Longrightarrow {\cI} =0 \q\text{\rm or} \q {\cJ} =0.
$$
On the other hand $\cA$ is said to be \textit{semiprime} if for
all ideals $\cI$ of $\cA$ we have
\[\cI \cI = 0 \Longrightarrow \cI = 0.\]
\end{definition}

The following lemma which is easily checked (see \cite[Exercise 1,
\S~8.2]{ZSSS}) tells us that the notion of $\cA$ being prime (or
semiprime) is independent of the base ring under which $\cA$ is
viewed as an algebra.

\begin{lemma}
\label{lem:primeinvariant} $\cA$ is prime (resp.~semiprime) as a
$B$-algebra if and only if $\cA$ is prime (resp.~semiprime) as a
$\bbZ$-algebra.
\end{lemma}

The following is proved  in \cite{EMO}.

\begin{lemma}
\label{lem:intdom} Assume $\cA$  is a prime algebra over $B$.
Then
\begin{itemize}
    \item[(i)] $C_B(\cA)$ is an
integral domain and $\cA$ is a torsion free $C_B(\cA)$-module.
    \item[(ii)] If we denote the quotient field
    of $C_B(\cA)$ by $\wt{C_B(\cA)}$, then
    $\cA\ot_{C_B(\cA)} \wt{C_B(\cA)}$ is a prime algebra over
    $\wt{C_B(\cA)}$.  Moreover, if $\cA$ is finitely generated as a module
    over its multiplication algebra
    $\Mult_B(\cA)$, then $\cA\ot_{C_B(\cA)} \wt{C_B(\cA)}$ is central
    over $\wt{C_B(\cA)}$.
\end{itemize}
\end{lemma}

\begin{proof} (i) is Theorem 1.1(a) of \cite{EMO},
whereas (ii) follows from Theorem 1.3(a) and (b) of \cite{EMO}.
\end{proof}

In  a later section of the paper we will investigate the
 type of an iterated loop algebra. In that section,
we will need the notion of central closure.

\begin{definition}
Let $\cA$ be a prime pfgc algebra over $B$.  Denote the quotient
field  of $C_B(\cA)$ by $\wt{C_B(\cA)},$ and form the
$\wt{C_B(\cA)}$-algebra
\begin{equation*}
\label{eqn:At} \tA := \cA\ot_{C_B(\cA)} \wt{C_B(\cA)}.
\end{equation*}
We call $\tA$ the  \emph{central closure of} $\cA$.
(This  is not apparently the same as the central closure defined
in \cite[\S~II]{EMO}.  Here we are following the terminology in,
for example, \cite[p.~154]{MZ}.) By Lemma~\ref{lem:intdom}(i),
$\cA$ is a torsion free $C_B(\cA)$-module, and so the map $a
\mapsto a\ot \one$ is an injection of $\cA$ into $\tA$ which we
regard as an identification. In this way \textit{$\cA$ is regarded
as a subalgebra of its central closure $\tA$}.
\end{definition}

We now  summarize the main facts that we will need about the
central closure:

\begin{proposition}
\label{prop:cquot1} Let $\cA$  be a prime pfgc algebra over $B$.
Then the central closure $\tA$
of $\cA$ is a prime pfgc algebra over $B$. Moreover, $\tA$ is
finite dimensional and central as an algebra over the
field~$\wt{C_B(\cA)}$.
\end{proposition}

\begin{proof}  $\tA$ is prime by Lemma \ref{lem:intdom}(ii).
Next, since $\cA$ is embedded as a subalgebra of $\tA$, we have
$\tA \ne 0$. Also, since $\cA$ is perfect, $\tA$ is perfect.
Furthermore, since $\cA$ is finitely generated as a
$C_B(\cA)$-module, $\tA$ is finitely generated as a
$\wt{C_B(\cA)}$-module and therefore also as a $C_B(\tA)$-module
(since $\lambda_{\tA}(\wt{C_B(\cA)}) \subset C_B(\tA))$.  Thus
$\tA$ is~pfgc.

We have just seen that $\tA$ is finite dimensional over
$\wt{C_B(\cA)}$.

Finally, since $\cA$ is pfgc, Proposition \ref{prop:fg}(ii) tells
us that $\cA$ is finitely generated as a $\Mult_B(\cA)$-module.
Thus $\tA$ is central over~$\wt{C_B(\cA)}$ by Lemma
\ref{lem:intdom}(ii).
\end{proof}

\section{Loop algebras}
\label{sec:loop}

\textbf{Assumptions and notation:} \textit{For the rest
of the article}, $k$ will denote a fixed base field.  Unless
indicated to the contrary, the term algebra will mean algebra over
$k$.  For the sake of brevity, if $\cA$ is an algebra (over $k$),
we will often write
\[\centA := C_k(\cA).\]

In this section we recall the definition of a loop  algebra and
derive some of its basic properties.

Throughout the section let $m$ be a positive integer and let
\[\Zmod{m} = \{\modi : i\in \bbZ\}\]
be the group of integers
modulo $m$,  where $\modi = i + m\bbZ \in \Zmod{m}$ for $i\in \bbZ$.
Let
\begin{equation*}
R=k[t^{\pm 1}] \andd S=k[z^{\pm 1}]
\end{equation*}
be the algebras of Laurent polynomials in the variables $t$ and
$z$ respectively, and we identify $R$ as a subalgebra of $S$ by
identifying
\[t = z^m.\]
Observe that $S$ is a free $R$-module of rank $m$ with basis
$\set{\one,z,\dots,z^{m-1}}$, and hence the ring extension $S/R$
is faithfully flat.

Recall that a $\Zmod{m}$-\textit{grading} of the algebra $\cA$ is an indexed family
$\Sigma = \set{\cA_\modi}_{\modi \in \Zmod{m}}$ of subspaces of $\cA$ so that
$\cA = \bigoplus_{\modi \in \Zmod{m}} \cA_\modi$ and
$\cA_\modi \cA_\modj \subset \cA_{\modi + \modj}$ for $\modi,\modj \in \Zmod{m}$.
The integer $m$ is called the \textit{modulus} of $\Sigma$.

\begin{definition}
\label{def:loop} Suppose that  $\cA$ is a $k$-algebra, and we are
given a  $\Zmod{m}$-grading $\Sigma$ of the algebra $\cA$:
\[
\cA= \bigoplus_{\modi \in \Zmod{m}} \cA_{\modi}.
\]
In $\cA\otk S$ we define
\[
L(\cA,\Sigma) := \bigoplus_{i\in \bbZ} \cA_{\modi}  \otk z^i =
(\cA_{\bar 0}  \otk R) \oplus (\cA_{\bar 1}  \otk zR) \oplus \dots
\oplus (\cA_{\overline{m-1}}  \otk z^{m-1} R) .
\]
Then $L(\cA,\Sigma)$ is an $R$-subalgebra of $\cA\otk S$ that we
call the \emph{loop algebra of $\Sigma$ based on~$\cA$}.  Since
$L(\cA,\Sigma)$ is an algebra over $R$,  $L(\cA,\Sigma)$ is also
an algebra over $k$.
\end{definition}

\begin{remark}
\label{rem:roleLaurent} If we wish to emphasize the role of the
variable $z$ in the construction of the loop algebra we
write $L(\cA,\Sigma)$ as $L(\cA,\Sigma,z)$.
\end{remark}

\begin{example}
If $m=1$, then $\Zmod{m} = \set{\bar 0}$,  $\cA = \cA_{\bar 0}$ and
$L(\cA,\Sigma) = \cA \otk S$ is called the \textit{untwisted} loop
algebra based on $\cA$.
\end{example}

\begin{remark}
\label{rem:primitive}
Suppose that $k$  contains a primitive $m^\thsup$
root of  unity $\zeta_m$. In that case we can choose to work with
finite order automorphisms of period $m$ rather than $\Zmod{m}$-gradings,
provided that we fix the choice of $\zeta_m$.

Indeed, suppose that $\cA$ is an algebra.
If $\sigma$ is an  algebra
automorphism of period $m$ of $\cA$, we may define a $\Zmod{m}$-grading
$\Sigma=\set{\cA_\modi}_{\modi \in \Zmod{m}}$ of $\cA$ by setting
\[
\cA_{\modi}=\set{x \in \cA \suchthat \sigma(x) = \zeta_m^i x },
\]
for $\modi \in \Zmod{m}$.  We refer to this grading $\Sigma$
as the \textit{grading determined by $\sigma$}.  It is clear
that any $\Zmod{m}$-grading is determined by a unique automorphism
$\sigma$ in this way.  If $\Sigma$ is the grading determined
by $\sigma$, we denote the
algebra $L(\cA,\Sigma)$ by $L(\cA,\sigma)$, or $L(\cA,\sigma,z)$
if we want to emphasize the role of $z$. The algebra
$L(\cA,\sigma)$ can alternately be defined as the subalgebra of
fixed points in $\cA\otk S$ of the automorphism $\sigma\ot
\eta_m^{-1}$, where $\eta_m\in Aut_k(S)$ is defined by $\eta_m(z)
= \zeta_m z$.
\end{remark}

\begin{remark}
\label{rem:history}
When $k = \mathbb C$, $\cA$ is a finite dimensional simple Lie
algebra over $k$ and $\sigma$ is a finite order automorphism of
$\cA$, the loop algebra $L(\cA,\sigma)$ was used by V.~Kac  in
\cite{K1} to give realizations of all affine Kac-Moody Lie
algebras and to classify finite order automorphisms of $\cA$. (See
\cite[Ch.~8]{K2} and \cite[Ch.~X, \S~5]{H} for more information
about this.)
\end{remark}

\textit{For the rest of the section, let  $\cA$ a $k$-algebra, let
$\Sigma$ be a grading of $\cA$ by $\Zmod{m}$, and let}
\[\cL = L(\cA,\Sigma).\]

We next describe a useful canonical form for elements of $\cA\otk
S$ in terms of elements of $L(\cA,\Sigma)$. For this purpose note
that $\cA\otk S$ is an $S$-module (with action denoted by
``$\cdot$'') and $L(\cA,\Sigma)$ is contained in $\cA\otk S$.
Thus we may write expressions like $\sum_{i=0}^{m-1} z^i \cdot x_i
\in \cA\otk S$ if $x_0,\dots,x_{m-1} \in L(\cA,\Sigma)$.

\begin{lemma}
\label{lem:canform1} Each element of $\cA\otk S$ can be written
uniquely in the form
\begin{equation}
\label{eq:canform1} \sum_{i=0}^{m-1} z^i \cdot x_i
\end{equation}
where $x_0,\dots,x_{m-1} \in \cL$.
\end{lemma}

\begin{proof}  This fact was proved using a Galois
cocycle argument in \cite[Theorem~3.6~(b)]{ABP2} in the case when
$k$ contains a primitive $m^\thsup$ root of unity.  We give a direct
proof here instead. Let $x\in \cA\otk S$. Then $x$ is the sum of
elements of the form $a\ot z^j$, where $j\in \bbZ$ and $a\in
\cA_{\bar \ell}$ for some $\ell\in \bbZ$.  But, if we write $j-\ell
= qm + i$, where $q\in \bbZ$ and $0\le i\le m-1$, then $a\ot z^j =
z^i\cdot(a\ot z^{j-i}) = z^i\cdot(a\ot z^{qm+\ell})$ and $a\ot
z^{qm+\ell}\in \cL$.  So $x$ can be expressed in the form
\eqref{eq:canform1}.  For uniqueness, suppose that
$\sum_{i=0}^{m-1} z^i \cdot x_i = 0$, where $x_0,\dots,x_{m-1} \in
\cL$.  Write $x_i = \sum_{j\in \bbZ} a_{ij}\ot z^j$, where
$a_{ij}\in \cA_{\modj}$ for all $j$ and only finitely many
$a_{ij}$ are non-zero. Then
\[\sum_{i=0}^{m-1} \sum_{j\in \bbZ} a_{ij}\ot z^{i+j} = 0.\]
For  $0\le \ell \le m-1$, the  $\cA_{\bar \ell}\otk S$-component
of the expression on the left above must be zero.  Thus we have
\[\sum_{i=0}^{m-1} \sum_{j  \equiv \ell} a_{ij}\ot z^{i+j} = 0\]
for $0\le \ell \le m-1$, where $\equiv$ denotes congruence modulo
$m$. The exponents $i+j$ appearing in this sum are all distinct
and so we have $a_{ij} = 0$ for all $i,j$ and hence $x_i = 0$ for
all  $i$.
\end{proof}

Next note that we have the canonical map $\xi  =\xi_{\cA,\Sigma} :
L(\cA,\Sigma)\ot_R S\to \cA\otk S$ defined by
\begin{equation*}
\xi\big(x\ot z^i\big) = z^i \cdot x
\end{equation*}
for  $x\in L(\cA,\Sigma)$, $i\in \bbZ$. As observed in
\cite[Theorem~3.6(b)]{ABP2}, Lemma \ref{lem:canform1} has the
following interpretation:

\begin{lemma}
\label{lem:form} The map $\xi_{\cA,\Sigma} : L(\cA,\Sigma)\ot_R
S\to \cA\otk S$ is an $S$-algebra isomorphism of
$L(\cA,\Sigma)\ot_R S$ onto $\cA\otk S$.
\end{lemma}

\begin{proof} Clearly $\xi$ is a homomorphism of
$S$-algebras.  Moreover, each element of $L(\cA,\Sigma)\ot_R S$
can be expressed in the form $\sum_{i=0}^{m-1}x_i \ot z^i$, where
$x_i\in L(\cA,\Sigma)$ for each~$i$, and so $\xi$ is bijective by
Lemma \ref{lem:canform1}.
\end{proof}

\begin{remark}
\label{rem:untwist1} Lemma \ref{lem:form}  tells us that after
base ring extension from $R$ to $S$ the loop algebra
$L(\cA,\Sigma)$ becomes isomorphic to the untwisted loop algebra
$\cA\otk S$. In other words, $L(\cA,\Sigma)$ is ``untwisted'' by
the extension $S/R$. This fact is of great importance in the study
of loop algebras since, among other things, it allows one to use
the  tools of Galois cohomology to study loop algebras \cite{ABP2, P}.
\end{remark}

We will need the following simple fact:

\begin{lemma}
\label{lem:perfect1}\
\begin{itemize}
\item[(i)] If $\cA\ne 0$, then $L(\cA,\Sigma)\ne 0$.
\item[(ii)] If $\cA$ is perfect, then $L(\cA,\Sigma)$ is perfect.
\end{itemize}
\end{lemma}
\begin{proof}
Statement (i) is clear and statement (ii) is easily checked (see
the argument  in \cite[Lemma~4.3]{ABP2}).
\end{proof}

We now examine the centroid  of $\cL = L(\cA,\Sigma)$.

First  note that since $\cL$  is an
$R$-algebra, $C_R(\cL)$ is naturally an $R$-algebra (see
Definition \ref{def:centroid}(iii)). So since $C_R(\cL) \subset
C(\cL)$, it follows that $C(\cL)$ is also an $R$-algebra.


Next  by Lemma~\ref{lem:graded} the centroid $\centA$ inherits a
$\Zmod{m}$-grading that we denote by~$C(\Sigma).$  Under this
grading we have
\[\centA =  \bigoplus_{\modi \in \Zmod{m}} \centA_{\modi},\]
where
\begin{equation}
\label{eq:gradeCdef} \centA_{\modi} = \set{\chi \in \centA
\suchthat \chi(\cA_{\modj})\subset \cA_{\modi + \modj} \text{ for
} \modj\in \Zmod{m}}.
\end{equation}

Now  let
\[
\psi  :=\psi_{\cA,\Sigma} : L(\centA,C(\Sigma))\to
C_R\big(L(\cA,\Sigma)\big)
\]
be the unique $k$-linear map so that
\[\big(\psi(\chi \ot z^i)\big)(a\ot z^j) = \chi(a)\ot z^{i+j}\]
for $i,j\in \bbZ$, $\chi\in \centA_\modi$, $a\in \cA_\modj$. It is
immediate from this definition that $\psi$ is a homomorphism of
$R$-algebras that we call canonical.

\begin{lemma}
\label{lem:centloop}  Assume $\cA$ is finitely generated as a
module over its multiplication algebra $\Mult_k(\cA)$. Then the
map $\psi_{\cA,\Sigma}: L\big(\centA,C(\Sigma)\big) \to
C_R\big(L(\cA,\Sigma)\big)$ is an $R$-algebra isomorphism.
\end{lemma}

\begin{proof}  Since the ring extension $S/R$ is faithfully flat, to
show that $\psi  $ is an $R$-module isomorphism it suffices to
show that  $\psi$ becomes an isomorphism of $S$-modules after the
base change   from $R$ to $S$. That this is so follows from the commutative
diagram
\[
\begin{array}{ccc}\\
L\big(\centA,C(\Sigma)\big)\ot_R S &\os{\psi  \ot \one}
\longrightarrow &C_R(\cL)\ot_R S\\
&&\downarrow  \nu_\cL\\
\xi_C\downarrow &&C_S\big(\cL\ot_R S\big)\\
&&\downarrow C_S(\xi)\\
\centA\otk S &\us{\nu_\cA}\longrightarrow &C_S(\cA\otk S)
\end{array}
\]
 in view of the fact that all vertical maps and the bottom row therein are
$S$-isomorphisms. In this diagram  $\xi_C =
\xi_{C({\cA}),C(\Sigma)}$ as in Lemma~\ref{lem:form}, $\nu_\cA =
\nu_{\cA,k,S}$ as in Lemma~\ref{lem:BaseChange}(a), $\nu_\cL =
\nu_{\cL,R,S}$ as in Lemma~\ref{lem:BaseChange}(b), and $C_S(\xi)$
is the isomorphism induced by the isomorphism $\xi =
\xi_{\cA,\Sigma} : \cL\ot_R S \to\cA\otk S$ (see Lemmas
\ref{lem:centinduce} and \ref{lem:form}).
\end{proof}

The following proposition tells us that the centroid of a loop
algebra based on a pfgc algebra $\cA$ is isomorphic to the loop
algebra of the centroid of $\cA$.

\begin{proposition}
\label{prop:centloop1}  Let $\cL = L(\cA,\Sigma)$ be a loop
algebra based on a pfgc algebra $\cA$. Then $C_R(\cL) = \centL$,
and the canonical map
\[\psi = \psi_{\cA,\Sigma} : L(\centA,C(\Sigma)) \to \centL\]
is an $R$-algebra isomorphism.
\end{proposition}

\begin{proof}
Since $\cL$ is perfect by Lemma \ref{lem:perfect1}(ii), it follows
that $C_R(\cL) = \centL$ by Lemma~\ref{lem:centroidinvariant}(ii).
Also since $\cA$ is pfgc, it follows from
Proposition~\ref{prop:fg}(ii) that $\cA$ is finitely generated as
a module over $\Mult_k(\cA)$. Thus, by Lemma~\ref{lem:centloop},
$\psi$ is an $R$-algebra isomorphism from $L(\centA,C(\Sigma))$
onto~$\centL$.
\end{proof}

Finally we want to show that a loop algebra based on a pfgc
algebra is pfgc. For this we will use the following:

\begin{lemma}
\label{lem:pfgc1} If $\cA$ is finitely generated as a
$\centA$-module then $L(\cA,\Sigma)$ is finitely generated as a
$C_R\big(L(\cA,\Sigma)\big)$-module.
\end{lemma}

\begin{proof}  Let $\{a_1,\dots,a_p\}$ be a set of homogeneous elements of $\cA$
that generates $\cA$ as a $\centA$-module.  Fix integers
$d_1,\dots,d_p$ so that $a_j \in \cA_{\overline{d_j}}.$ Let
$\mathcal M$ be the $C_R(\cL)$-submodule of $\cL$ generated by the
elements $a_k\ot z^{d_k}.$ Since $S/R$ is flat we may identify
$\mathcal M \ot_R S$ as an $S$-submodule of $\cL \ot_R S$, and
since $S/R$ is faithfully flat it is sufficient to show that
$\mathcal M \ot_R S = \cL \ot_R S$ \cite[Ch.~I, \S~3.1,
Proposition~2]{B:CA}. We do this by showing that $\xi(\mathcal M
\ot_R S) = \xi(\cL \ot_R S)$, where $\xi = \xi_{\cA,\Sigma}
:\cL\ot_R S \longrightarrow \cA\otk S$ is the $S$-algebra
isomorphism from Lemma~\ref{lem:form}.

Suppose that $i,j\in \bbZ$, $\chi\in \centA_\modi$ and $1\le \ell
\le p$.  Then $\psi(\chi\ot z^i)$ is an element of $C_R(\cL)$,
where $\psi = \psi_{\cA,\Sigma}$. So $\big(\psi(\chi\ot
z^i)\big)(a_\ell\ot z^{d_\ell})\in \mathcal M$. But under $\xi$
we have
\begin{equation*}
\Big(\big(\psi(\chi\ot z^i)\big)(a_\ell\ot z^{d_\ell})\Big)\ot
z^j \mapsto \chi(a_\ell)\ot z^{d_\ell+i+j}.
\end{equation*}
Since $\{a_1,\dots,a_p\}$ generates $\cA$ as a $\centA$-module, it
follows $\xi(\mathcal M\ot_R S) =\xi(\cL\ot_R S)$ as needed.
\end{proof}

\begin{proposition}
\label{prop:pfgc1}  Let $\cL = L(\cA,\Sigma)$ be a loop algebra
based on a pfgc algebra $\cA$. Then $\cL$ is a pfgc algebra.
\end{proposition}

\begin{proof} $\cL \ne (0)$ and $\cL$ is perfect
by Lemma~\ref{lem:perfect1}.  So P0 and P1 hold (see Definition
\ref{def:pfgc}). By Lemma~\ref{lem:pfgc1}, $\cL$ is finitely
generated as a $C_R(\cL)$-module. But by
Proposition~\ref{prop:centloop1}, we have $C_R(\cL) = \centL$.
Thus $\cL$ is finitely generated as a $\centL$-module and so P2
holds.  Hence $\cL$ is pfgc.
\end{proof}

\section{Iterated loop algebras}
\label{sec:iterated}

In  this section we define iterated loop algebras and prove some
of their basic properties.

{\bf Notation:} \textit{For the  rest of this article, we fix some notation.}
Let $n$ be a positive integer. Let
$z_1,\dots,z_n$ be a sequence of algebraically independent variables
over~$k$. For $0\le p \le n$, let
\[\Stp := k[z_1^{\pm 1},\dots,z_p^{\pm 1}]\]
be the algebra of Laurent polynomials in the variables
$z_1,\dots,z_p$ over~$k$.  (So $S^{\ot 0} = k$.) We identify
$S^{\ot p} \ot S^{\ot q} = S^{\ot (p+q)}$ in the natural fashion
when $0\le p,q\le n$ and $p+q\le n$. We also fix a
sequence $m_1,\dots,m_n$ of positive integers, and we set
\[I_p := \set{(i_1,\dots,i_p) \in \bbZ^p \suchthat
0\le i_j \le m_j-1 \text{ for all } j},\]
for  $1\le p \le n$.

\begin{definition}
\label{def:iteratedloop} Suppose that $\cA$ is an algebra over
$k$. An algebra $\cL$ over $k$ is called an \emph{$n$-step loop
algebra}  or an \emph{iterated loop algebra} based on $\cA$ if
there exists a sequence $\cL_0, \cL_1,\dots,\cL_n$ of algebras so
that $\cL_0 = \cA$, $\cL_n = \cL$ and
\[\cL_{p} = L(\cL_{p-1},\Sigma_{p},z_{p}),\]
for $1\le p \le n$, where $\Sigma_{p}$ is a  $\Zmod{m_p}$-grading
of $\cL_{p-1}$. (See Remark \ref{rem:roleLaurent} for the notation
used here.) In that case we write
\[\cL = L(\cA,\Sigma_1,\dots,\Sigma_n)\]
(suppressing in the notation the role of the  variables
$z_1,\dots,z_n$).
\end{definition}

\begin{remark}
\label{rem:iteratedloop}  Suppose that $\cL$ is an $n$-step loop
algebra based on $\cA$ and we have the notation from Definition
\ref{def:iteratedloop}.

(i) For $1 \le p \le n$, $\cL_p =
L(\cA,\Sigma_1,\dots,\Sigma_p)$ is a $p$-step loop algebra
based on $\cA$.

(ii)  Observe that $\cL_{p} \subset \cL_{p-1}\otk k[z_{p}^{\pm
1}]$ for $1\le p \le n$.  Thus
\[\cL_p \subset (\dots((\cA \otk k[z_1^{\pm 1}])\otk k[z_2^{\pm 1}])\dots)\otk k[z_p^{\pm 1}]
= \cA\ot k[z_1^{\pm 1},\dots,z_p^{\pm 1}],\] for $0\le p \le n$,
where the last equality is the natural identification using the
associativity of the tensor product and the identification
$k[z_1^{\pm 1}]\otk \dots \otk k[z_p^{\pm 1}] = \Stp$.
Consequently, \textit{$\cL_p$ is a subalgebra of $\cA\otk \Stp$
for $0\le p \le n$, and in particular $\cL$ is a subalgebra of
$\cA\otk \Stn$}.

(iii) Suppose that $k$ contains a primitive $m_i^\thsup$ root of
unity $\zeta_{m_i}$ (which we fix) for $i=1,\dots,n$.  Then for $p=
1,\dots,n$, the grading $\Sigma_p$ of $\cL_{p-1}$ is determined by a
unique automorphism $\sigma_p$ of  $\cL_{p-1}$ of period
$m_p$. We then denote the algebra $L(\cA,\Sigma_1,\dots,\Sigma_n)$
by $L(\cA,\sigma_1,\dots,\sigma_n)$.
\end{remark}

\begin{example}
If $m_1 = \dots m_n=1$ then $L(\cA,\Sigma_1,\dots,\Sigma_n) = \cA
\otk \Stn$ is called the \textit{untwisted} $n$-step loop algebra
based on $\cA$.
\end{example}

\begin{example}[Multiloop algebras]
\label{ex:multiloop}
Suppose  that $k$ contains a primitive ${m_i}^\thsup$ root of unity
$\zeta_{m_i}$ for $1\le i \le n$.  Let $\cA$ be an algebra, and let
$\sigma_1,\dots,\sigma_n$ be commuting finite order automorphisms of
$\cA$ with periods $m_1,\dots,m_n$ respectively.  Let
\[\cA_{\modi_1,\dots,\modi_n} = \set{x\in \cA \suchthat
\sigma_j x = \zeta_{m_j}^{i_j} x \text{ for } 1\le j\le n}\]
for $(i_1,\dots,i_n)\in \bbZ^n$, where $\modi_j := i_j + m_j\bbZ \in \Zmod{m_j}$
for $1\le j\le n$.  Then
\[\cA = \bigoplus_{(i_1,\dots,i_n)\in I_n} \cA_{\modi_1,\dots,\modi_n},\]
and we set
\[\Mloop(\cA,\sigma_1,\dots,\sigma_n) := \bigoplus_{(i_1,\dots,i_n)\in \bbZ^n}
\cA_{\modi_1,\dots,\modi_n}\otk z_1^{i_1}\dots z_n^{i_n}\]
in $\cA\otk \Stn$.  Then $\Mloop(\cA,\sigma_1,\dots,\sigma_n)$ is a subalgebra of
$\cA\otk \Stn$ that we call the $n$-step \textit{multiloop algebra} of $\sigma_1,\dots,\sigma_n$
based on $\cA$.

Now the multiloop algebra $\cL = \Mloop(\cA,\sigma_1,\dots,\sigma_n)$ can be interpreted
as an iterated loop algebra.  To see this, let $\cL_0 = \cA$ and let
$\cL_p = \Mloop(\cA,\sigma_1,\dots,\sigma_p)$ for $1\le p\le n$.
Then by definition we have $\cL_0 = \cA$ and $\cL_n = \cL$.  Also,
for $1\le p \le n$, we may define a $\Zmod{m_{p}}$-grading $\Sigma_{p}$
on
$\cL_{p-1}$ by setting
\[(\cL_{p-1})_{\modi_{p}} =  \bigoplus_{(i_1,\dots,i_{p-1})\in \bbZ^{p-1}}
\cA_{\modi_1,\dots,\modi_{p}}\otk z_1^{i_1}\dots z_{p-1}^{i_{p-1}}\]
for $\modi_{p}\in \Zmod{m_{p}}$, in which case
it is then clear that $L(\cL_{p-1},\Sigma_{p},z_{p}) = \cL_{p}$.
So $\cL = L(\cA,\Sigma_1,\dots,\Sigma_n)$.
\end{example}

We have just seen in Example
\ref{ex:multiloop} that any multiloop algebra is an iterated loop algebra.
However we will see later in Example \ref{ex:hermitian} that there are iterated loop algebras
$\cA$ that are not multiloop algebras.

For the rest of the section \textit{we assume that
\[\cL = L(\cA,\Sigma_1,\dots,\Sigma_n)\]
is an $n$-step loop algebra based on an algebra $\cA$ over $k$},
and we use the notation $\cL_0,\dots,\cL_n$ of Definition
\ref{def:iteratedloop}.

Our  first theorem describes some important basic algebraic
properties that are inherited by a loop algebra from its base. In
the last part of this theorem we will see how the  Krull dimension
of a loop algebra depends on the Krull dimension of its base. Here
and subsequently we use
\[\Dim \mathcal C\]
to denote the \textit{Krull dimension} of a unital commutative
associative $k$-algebra $\mathcal C$ (when regarded as a ring).
Note that if $\mathcal C$ is finitely generated as a $k$-algebra
then $\Dim \mathcal C$ is finite  \cite[p.~52]{Ku}.

\begin{theorem} Let $\cL = L(\cA,\Sigma_1,\dots,\Sigma_n)$.
\label{thm:basic}
\begin{itemize}
\item[(i)] If $\cA\ne 0$ then $\cL\ne 0$.
\item[(ii)]
If $\cA$ is perfect then $\cL$ is perfect.
\item[(iii)]
If $\cA$ is pfgc then $\cL$ is pfgc.
\item[(iv)] If $\cA$ is prime then $\cL$ is prime.
\item[(v)] If $\cA$ is unital then $\cL$ is a unital subalgebra of $\cA\ot \Stn$.
\item[(vi)] If
$\cA$ is commutative then $\cL$ is commutative.
\item[(vii)]
If $\cA$ is associative then $\cL$ is associative.
\item[(viii)]  If $\cA$ is an integral domain
then $\cL$ is an integral domain.
\item[(ix)] If $\cA$ is unital and finitely generated as a $k$-algebra
then $\cL$ is unital and finitely generated as a $k$-algebra.
\item[(x)]
If $\cA$ is unital, commutative, associative and finitely
generated as a $k$-algebra, then $\cL$ is unital, commutative,
associative and finitely generated as a $k$-algebra and
\begin{equation}
\label{eq:Krull} \Dim \cL = \Dim \cA +n.
\end{equation}
\end{itemize}
\end{theorem}

\begin{proof}  Since $\cL_{p+1}$ is a loop algebra based on $\cL_p$
for $0\le p\le n-1$, we can assume in the proof of each of these
statements that $n=1$.  So we may use the notation of
Section~\ref{sec:loop}:
\[m=m_1,\ z=z_1,\ \Sigma = \Sigma_1,\ \cL = L(\cA,\Sigma,z),\
S = S^{\ot 1} = k[z^{\pm 1}] \text{ and } R = k[z^{\pm m}].\]

Now (i) and (ii) follow from Lemma \ref{lem:perfect1}. (iii)
follows from Proposition \ref{prop:pfgc1}. (v) follows from that
fact that $\one_\cA\in \cA_{\bar 0}$, since then $\one_\cA \ot
1_S\in \cL$. (vi), (vii) and (viii) follow from the fact that
$\cL$ is a subalgebra of $\cA\otk S \simeq_k \cA[z^{\pm 1}]$. So
we only need to prove (iv), (ix) and (x).

(iv) We  show first that $\cA\otk S$ is prime. For this let $\cI$
and $\cJ$ be ideals of the $k$-algebra $\cA\otk S$ such that
$\cI\cJ = 0$.  For $m\in \bbZ$, let
\[\cI_m = \set{a\in \cA \suchthat
\exists \ a_i\in \cA \text{ for } i\ge m \text{ with } a_m = a
\text{ and } \textstyle{\sum_{i\ge m}} a_i\ot z^i \in \cI  },\] in
which case $\cI_m$ is an ideal of $\cA$. Similarly, using $\cJ$
instead of $\cI$, we define an ideal $\cJ_n$ of $\cA$ for $n\in
\bbZ$.  Furthermore, since $\cI \cJ = 0$, we have $\cI_m\cJ_n = 0$
for $m,n\in \bbZ$. Now suppose that $\cI \ne 0$.  Then $\cI_m \ne 0$
for some $m\in Z$.  Thus, since $\cA$ is prime, we have $\cJ_n =
0$ for all $n\in \bbZ$ and so $\cJ = 0$. Therefore $\cA\otk S$ is
prime.

But  $\cL\ot_R S\simeq_S \cA\otk S$ by Lemma \ref{lem:form}. Hence
$\cL\ot_R S$ is a prime algebra. To prove that $\cL$ is prime (as
a $k$-algebra), it is enough to show that $\cL$ is a prime
$R$-algebra (by Lemma \ref{lem:primeinvariant}). For this let
$\cI$ and $\cJ$ be ideals of the $R$-algebra $\cL$ such that
$\cI\cJ=0$.  Since $S/R$ is flat, we can identify $\cI\ot_R S$ and
$\cJ\ot_R S$ as ideals of the $S$-algebra $\cL\ot_R S$.
Furthermore, we have $(\cI\ot_R S)(\cJ\ot_R S)=0.$  Since
$\cL\ot_R S$ is prime, either $\cI\ot_R S$ or $\cJ\ot_R S$ is $0$.
Therefore $\cI=0$ or $\cJ=0$ by the faithful flatness of $S/R$.

(ix) $\cL$ is unital by (v). Let $\{a_1,\dots,a_p\}$ be a set of
homogeneous elements of $\cA$ that generates $\cA$ as a
$k$-algebra, and fix integers  $d_1,\dots,d_p$ so that $a_j \in
\cA_{\bar d_j}.$ One easily checks that the elements
$a_1\ot z^{d_1},\dots,a_p\ot z^{d_p}$ together with the elements
$\one_\cA\ot z^m$ and $\one_\cA\ot z^{-m}$ generate $\cL$ as a
$k$-algebra.

(x) $\cL$ is  unital, commutative, associative and finitely
generated as a $k$-algebra by (v), (vi), (vii) and (ix), and so
the Krull dimensions of both $\cA$ and $\cL$ are finite. Now
recall  that $\cL$ is a subalgebra of $\cA\otk S$ and, by Lemma
\ref{lem:canform1}, each element of $\cA\otk S$ can be written
uniquely in the form
\[\sum_{i=0}^{m-1} z^i \cdot x_i =  \sum_{i=0}^{m-1}x_i(\one_\cA \ot z^i) \]
where $x_0,\dots,x_{m-1} \in \cL$.  Thus $\cA\otk S$ is a free
$\cL$-module of rank $m$, and  so in particular $\cA\otk S$
is a finitely generated $\cL$-module.  Hence $\cA \otk S/\cL$ is
an integral ring extension and so by \cite[Corollary II.2.13]{Ku},
\[
\Dim \cL = \Dim \big(\cA\otk S\big).
\]
On the other hand since both $\cA$ and $S$ are finitely generated
$k$-algebras
\[
\Dim \big(\cA\otk S\big) = \Dim \cA + \Dim S
\]
\cite[Corollary II.3.9]{Ku}. Since $\Dim S =1$,
we obtain  $\Dim \cL
=\Dim \cA +1.$
\end{proof}

\begin{remark}
\label{rem:primepfgc} It follows in particular from Theorem
\ref{thm:basic} that any $n$-step loop algebra based on a prime
pfgc algebra is a prime pfgc algebra. The corresponding statement
is not true for simple pfgc algebras. For example, an untwisted
pfgc algebra $\cA\otk \Stn$ is never simple (since $\Stn$ is not
simple).  This is the reason why prime pfgc algebras are natural
algebras to consider when studying loop algebras, even if one's
main interest is in the case when the base algebras are simple.
\end{remark}

We  conclude this section with a generalization to iterated loop
algebras of the canonical form described in Lemma
\ref{lem:canform1}. If $1\le p \le n$,  we use the usual convenient
notation
\[z^{\bi} = z_1^{i_1}\dots z_p^{i_p}\]
for $\bi = (i_1,\dots,i_p) \in \bbZ^p$.
Note that $\cA\otk \Stp$ is an $\Stp$-module (with action
denoted by ``$\cdot$''), and $\cL_p$ is contained in $\cA\otk
\Stp$.  Thus, we can write expressions like $\sum_{\bi\in I_p}
z^\bi\cdot x_\bi\in \cA\otk \Stp$, where $x_\bi\in \cL_p$ for all
$\bi\in I_p$.

\begin{lemma}
\label{lem:canform} If  $1\le p\le n$, each element in  $\cA\otk
\Stp$ can be expressed uniquely in the form
\begin{equation}
\label{eq:canform} \sum_{\bi\in I_p} z^\bi\cdot x_\bi,
\end{equation}
where $x_\bi\in \cL_p$ for all $\bi$.
\end{lemma}

\begin{proof}  We argue by induction on $p$.  When $p=1$, the statement
follows from Lemma \ref{lem:canform1}. So we suppose that the
statement is true for $p$, where $1\le p\le n-1$.

Let $x\in \cA\ot \Stpplus$.  To show that $x$ can be expressed in
the form \eqref{eq:canform}, note first that $x$ is a sum of
elements of the form $x'\ot z_{p+1}^j$, where $x'\in \cA\otk\Stp$
and $j\in \bbZ$. But by the induction hypothesis, $x'$ is the sum of
elements of the form $z^\bi\cdot x''$, where $\bi\in I_p$ and
$x''\in \cL_p$. Thus $x$ is the sum of elements of the form
\[(z^\bi\cdot x'')\ot z_{p+1}^j = z^\bi\cdot (x''\ot z_{p+1}^j). \]
But $x''\ot z_{p+1}^j\in \cL_p\otk k[z_{p+1}^{\pm 1}]$, and so, by
Lemma \ref{lem:canform1}, $x''\ot z_{p+1}^j$ is the sum of
elements of the form $z_{p+1}^\ell \cdot x'''$, where $0\le \ell
\le m_{p+1}-1$ and $x'''\in \cL_{p+1}$.  Thus $x$ is the sum of
elements of the form
\[z^\bi\cdot(z_{p+1}^\ell \cdot x''') = (z^\bi z_{p+1}^\ell)\cdot x''' \]
as desired.

For uniqueness, suppose that $\sum_{\bj\in I_{p+1}} z^\bj\cdot
x_\bj=0$, where $x_\bj\in \cL_{p+1}$ for each $\bj\in I_{p+1}$.
Then
\[\sum_{\bi\in I_p} \sum_{\ell=0}^{m_{p+1}-1 }(z^\bi z_{p+1}^\ell)\cdot x_{\bi,\ell} = 0,\]
where, if $\bi=(i_1,\dots,i_p)\in I_p$ and $0\le \ell\le
m_{p+1}-1$, we are using the notation $x_{\bi,\ell} :=
x_{(i_1,\dots,i_p,\ell)}\in \cL_{p+1}$. So we have
\[\sum_{\bi\in I_p} z^\bi \cdot
\bigg(\sum_{\ell=0}^{m_{p+1}-1 } z_{p+1}^\ell\cdot
x_{\bi,\ell}\bigg) = 0.\] But for $\bi\in I_p$, the element
$\sum_{\ell=0}^{m_{p+1}-1 } z_{p+1}^\ell\cdot x_{\bi,\ell}$ is in
$\cL_p\otk k[z_{p+1}^{\pm 1}]$ and therefore we can write
\[\sum_{\ell=0}^{m_{p+1}-1 } z_{p+1}^\ell\cdot x_{\bi,\ell} =
\sum_{j\in \bbZ} y_{\bi,j}\ot z_{p+1}^j,\] where each $y_{\bi,j}$ is
in $\cL_p$ and only finitely many of these elements  are nonzero.
Then $\sum_{\bi\in I_p} z^\bi \cdot \big(\sum_{j\in \bbZ}
y_{\bi,j}\ot z_{p+1}^j\big) = 0$, and so
\[ \sum_{j\in \bbZ}
\bigg(\sum_{\bi\in I_p} z^\bi \cdot y_{\bi,j}\bigg) \ot z_{p+1}^j
= 0.\] Hence $ \sum_{\bi\in I_p} z^\bi \cdot y_{\bi,j}=0$ for each
$j$ and so by the induction hypothesis $y_{\bi,j}=0$ for all
$\bi\in I_p$ and~$j\in \bbZ$. So $\sum_{\ell=0}^{m_{p+1}-1 }
z_{p+1}^\ell\cdot x_{\bi,\ell} = 0$ for all $\bi\in I_p$, and
hence, by Lemma \ref{lem:canform1}, $x_{\bi,\ell} = 0$ for all
$\bi\in I_p$ and $0\le \ell \le m_{p+1}-1$.
\end{proof}

If $\cA$ is unital and associative, then $\cL =
L(\cA,\Sigma_1,\dots,\Sigma_n)$ is a unital associative
subalgebra of $\cA \otk \Stn$ and hence $\cA \otk \Stn$ is an
$\cL$-module (with action denoted by~``$\cdot$'').

\begin{corollary}
\label{cor:freemod} Suppose that $\cL =
L(\cA,\Sigma_1,\dots,\Sigma_n)$ where $\cA$ is unital and
associative. Then $\cA\otk \Stn$ is a free $\cL$-module of rank
$m_1\dots m_n$ with basis
\[\set{\one_\cA\ot z^\bi}_{\bi\in I_n}.\]
\end{corollary}

\begin{proof}
This follows from Lemma \ref{lem:canform} (with $p=n$) and the
observation that
\[z^\bi \cdot x = x\cdot (\one_\cA\ot z^\bi)\]
for $x\in \cL$ and $\bi\in \bbZ^n$. (On the left of this equation
``$\cdot$'' denotes the action of $\Stn$ on $\cA\otk \Stn$,
whereas on the right ``$\cdot$'' denotes the action of $\cL$ on
$\cA\otk \Stn$.)
\end{proof}

\section{The centroid of an iterated loop algebra}
\label{sec:centroiditerated}

In this section,  we give an explicit  description of the centroid
of an $n$-step loop algebra based on a pfgc algebra $\cA$ as an
$n$-step loop algebra based on $\centA$.

Throughout  the section \textit{we assume that $\cL =
L(\cA,\Sigma_1,\dots,\Sigma_n)$ is an $n$-step loop algebra
based on a algebra $\cA$ over $k$}.  So we have algebras
$\cL_0,\dots,\cL_n$ so that $\cL_0 = \cA$, $\cL_n = \cL$ and
\[\cL_{p+1} = L(\cL_p,\Sigma_{p+1},z_{p+1})\]
for $0\le p \le n-1$.  As we observed in Remark
\ref{rem:iteratedloop}, $\cL_p$ is a subalgebra of $\cA \otk \Stp$
for $0\le p\le n$.

We next introduce some notation.

First let $0\le p\le n$.  Then $\centA\otk\Stp$ is a unital
associative algebra and \textit{$\cA\otk \Stp$ is a
$\centA\otk\Stp$-module} under the action ``$\cdot$'' defined by
\begin{equation*}
(\chi \ot z^\bi)\cdot (a \ot z^\bj) = \chi(a) \ot z^{\bi + \bj}.
\end{equation*}
We let $\bC(\cL_p)$ denote the stabilizer of $\cL_p$ in
$\centA\otk\Stp$ under this action.  That is we let
\[\bC(\cL_p) := \set{u\in \centA\otk\Stp \suchthat u\cdot \cL_p \subset \cL_p}.
\]
Then $\bC(\cL_p)$ is a unital subalgebra of $\centA\otk\Stp$ and
$\cL_p$ is a $\bC(\cL_p)$-module. (For convenience, our notation
suppresses the fact that $\bC(\cL_p)$ depends on $\cA$,
$\Sigma_1,\dots,\Sigma_p$ and not just on the loop algebra
$\cL_p$.)

Next suppose that $0\le p\le n-1$.  Then $\Sigma_{p+1}$ is a
$\Zmod{m_{p+1}}$-grading of the algebra $\cL_p$ which we write as
\[ \cL_p = \bigoplus_{\modi \in \Zmod{m_{p+1}}} (\cL_p)_\modi.\]
We set
\begin{equation}
\label{eq:gradebCdef} \bC(\cL_p)_\modi := \set{u\in \bC(\cL_p)
\suchthat u\cdot (\cL_p)_{\modj}
 \subset (\cL_p)_{\modi + \modj} \text{ for all } \modj\in \Zmod{m_{p+1}}}
\end{equation}
for $\modi \in \Zmod{m_{p+1}}$.  We denote the collection
$\{\bC(\cL_p)_{\modi}\}_{\modi \in \Zmod{m_{p+1}}}$ by $\bC(\Sigma_{p+1})$.
We will see in Lemma \ref{lem:bCp}(ii) below that $\bC(\Sigma_{p+1})$
is a $\Zmod{m_{p+1}}$-grading of $\bC(\cL_p)$.

Finally for $0\le p \le n$ we define $\gamma_p : \bC(\cL_p) \to
C(\cL_p)$ by
\[\gamma_p(u)(x) = u\cdot x\]
for $u\in \bC(\cL_p)$, $x\in \cL_p$, in which case $\gamma_p$ is
an algebra homomorphism. Note in particular that $\bC(\cL_0) =
\centA$ and $\gamma_0$ is the identity map.

\begin{lemma}
\label{lem:bCp} Suppose that $\cL =
L(\cA,\Sigma_1,\dots,\Sigma_n)$, where $\cA$ is a pfgc
algebra.
\begin{itemize}
    \item[(i)] If $0\le p\le n$ then $\gamma_p : \bC(\cL_p) \to C(\cL_p)$ is an
    isomorphism of $k$-algebras.
    \item[(ii)] If $0\le p\le n-1$ then $\bC(\Sigma_{p+1})$ is a $\Zmod{m_{p+1}}$-grading
    of the algebra $\bC(\cL_p)$ and the map $\gamma_p$ is an isomorphism of graded
    algebras.
    \item[(iii)] If $0\le p\le n-1$ then
        \begin{equation}
        \label{eq:barCiterate}
        \bC(\cL_{p+1}) = L(\bC(\cL_p),\bC(\Sigma_{p+1}),z_{p+1}).
        \end{equation}
\end{itemize}
\end{lemma}

\begin{proof} (i) We first  show that
$\gamma_p$ is injective for $0\le p\le n$.  To see this, suppose
that $u\in \ker(\gamma_p)$. Then $u\cdot x = 0$ for all $x\in
\cL_p$, and so (since $\cL_p$ spans $\cA\otk \Stp$ over $\Stp$ by
Lemma \ref{lem:canform}) we have $u\cdot x = 0$ for all $x\in
\cA\otk \Stp$. This implies that $u=0$.

Next we prove the bijectively of $\gamma_p$ for $0\le p\le n$ by
induction on $p$.  This is clear if $p=0$ since $\gamma_0$ is the
identity map. So we suppose that $0\le p \le n-1$ and that
$\gamma_p$ is bijective. It is clear from this bijectivity and
from the definitions of $\bC(\cL_{p+1})$ and
$C(\cL_{p+1})$ (see
\eqref{eq:gradebCdef} and \eqref{eq:gradeCdef}) that
\[\gamma_p(\bC(\cL_{p+1})_\modi) = C(\cL_{p+1})_\modi\]
for $\modi\in \Zmod{m_{p+1}}$. Hence $\bC(\Sigma_{p+1})$ is a grading of the
algebra $\bC(\cL_p)$ and $\gamma_p : \bC(\cL_p)\to C(\cL_p)$ is a
graded isomorphism. So $\gamma_p$ induces an algebra isomorphism
\[ L(\gamma_p) : L(\bC(\cL_p),\bC(\Sigma_{p+1}), z_{p+1}) \to
L(C(\cL_p),C(\Sigma_{p+1}),z_{p+1}).\] Consequently we have the
composite algebra isomorphism
\begin{equation}
\label{eq:comp} L(\bC(\cL_p),\bC(\Sigma_{p+1}), z_{p+1})
\os{L(\gamma_p)} \longrightarrow L(C(\cL_p),C(\Sigma_{p+1}),z_{p+1})
\os{\psi_{\cL_p,\Sigma_{p+1}}} \longrightarrow C(\cL_{p+1}),
\end{equation}
where $\psi_{\cL_p,\Sigma_{p+1}}$ is the isomorphism of  Proposition
\ref{prop:centloop1}.  (Note that Proposition \ref{prop:centloop1}
can be applied since $\cL_p$ is a pfgc algebra by Theorem
\ref{thm:basic}(iii).) But
\[L(\bC(\cL_p),\bC(\Sigma_{p+1}), z_{p+1}) \subset \bC(\cL_{p+1})\]
and one easily checks that the restriction
\begin{equation}
\label{eq:restrict}
\gamma_{p+1} |_{L(\bC(\cL_p),\bC(\Sigma_{p+1}),z_{p+1})}
: L(\bC(\cL_p),\bC(\Sigma_{p+1}),z_{p+1})
\to C(\cL_{p+1})
\end{equation}
equals the composite map \eqref{eq:comp}. Hence the restriction
\eqref{eq:restrict} of $\gamma_{p+1}$ is bijective. Thus, since
$\gamma_{p+1}$ itself is injective, it follows that
\[L(\bC(\cL_p),\bC(\Sigma_{p+1}), z_{p+1}) = \bC(\cL_{p+1})\]
and $\gamma_{p+1}$ is bijective.  So we have proved (i).

(ii) and (iii):  These were proved in the argument for (i).
\end{proof}

Since $\cL = \cL_n$, we write $\bC(\cL) = \bC(\cL_n)$ and so
\[\bC(\cL) := \set{u\in \centA\otk S^{\ot n}
\suchthat u\cdot \cL \subset \cL}.
\]
Then $\bC(\cL)$ is a unital subalgebra of $\centA\otk S^{\ot n}$,
and $\cA \ot_k S^{\ot n}$ is a $\bC(\cL)$-module. We also write
$\gamma_\cL  = \gamma_n$. Thus $\gamma_\cL  : \bC(\cL) \to \centL$
is the $k$-algebra homomorphism (said to be canonical) defined by
\[\gamma_\cL(u)(x) = u\cdot x\]
for $u\in \bC(\cL)$, $x\in \cL$.

Using Lemma \ref{lem:bCp} we can now give an explicit description
of the centroid of an $n$-step loop algebra as an $n$-step loop
algebra.

\begin{theorem}
\label{thm:centloopn} Suppose that $\cL =
L(\cA,\Sigma_1,\dots,\Sigma_n)$ is an $n$-step loop algebra
based on a pfgc algebra $\cA$. Then the canonical map $\gamma_\cL
: \bC(\cL) \to \centL$ is an algebra isomorphism and we have
\begin{equation}
\label{eq:barCloop} \bC(\cL) =
L(\centA,\bC(\Sigma_1),\dots,\bC(\Sigma_n)).
\end{equation}
\end{theorem}

\begin{proof}
$\gamma_\cL$ is an isomorphism by Lemma \ref{lem:bCp}(i).
Moreover \eqref{eq:barCloop} follows by repeated application of
\eqref{eq:barCiterate}.
\end{proof}

\goodbreak

\begin{corollary}
\label{cor:Krullcent} Suppose that $\cL$ is an $n$-step loop
algebra based on a pfgc algebra $\cA$. Then
\begin{itemize}
\item[(i)] $\centA$ and $C(\cL)$ are nonzero unital commutative
associative algebras over~$k$.
\item[(ii)]  If $\centA$ is an integral domain, then
$C(\cL)$ is an integral domain.
\item[(iii)]  If $\centA$ is finitely generated as an algebra
over $k$, then $C(\cL)$ is finitely generated as an algebra over
$k$ and $\Dim C(\cL) = \Dim \centA +n$.
\end{itemize}
\end{corollary}

\begin{proof}
(i)  Since $\cA$ is pfgc, we know that $\cL$ is pfgc by Theorem
\ref{thm:basic}(iii). Hence $\centA$ and $C(\cL)$ are nonzero
unital commutative associative algebras by
Proposition~\ref{prop:fg}(i).

(ii) and (iii):  We know by Theorem \ref{thm:centloopn} that
$C(\cL)$ is isomorphic to an $n$-step loop algebra based on
$\centA$. Thus (ii) and (iii) follow from Theorem \ref{thm:basic}
(viii) and (x) respectively.
\end{proof}

If $\cA$ is a finite dimensional central simple algebra over $k$,
then $\centA = k$ and $\cA$ is a pfgc algebra.  Hence we have the
following consequence of Corollary \ref{cor:Krullcent}:

\begin{corollary}
\label{cor:centsimplebasea} Suppose that $\cL$ is an $n$-step loop
algebra based on a finite dimensional central simple algebra
$\cA$ over $k$. Then $C(\cL)$ is an integral domain, $C(\cL)$ is finitely
generated as an algebra over $k$, and $\Dim C(\cL) = n$.
Consequently, if $\cL'$ is an $n'$-step loop algebra based on a
finite dimensional central simple algebra $\cA'$ over $k$, then
\[\cL \simeq_k \cL' \implies n = n'.\]
\end{corollary}

\begin{remark}
\label{rem:centcalc}
Suppose that $\cL$ is an $n$-step loop algebra based on
a finite dimensional central simple algebra $\cA$  over $k$.  Then  $\centA
\otk \Stn = k\otk \Stn = \Stn$ and so
\[C(\cL) \os{\gamma_\cL}\simeq_k \bC(\cL) = \set{u\in \Stn : u\cdot \cL \subset \cL}.\]
This fact can be used to explicitly compute $C(\cL)$ in examples.
\end{remark}

\begin{corollary}
\label{cor:centmulti}
Suppose  that $\cL = \Mloop(\cA,\sigma_1,\dots,\sigma_n)$ is a multiloop
algebra based on a finite dimensional central
simple algebra $\cA$  over $k$, where $\sigma_1,\dots,\sigma_n$
are commuting finite order automorphisms of $\cA$ with periods
$m_1,\dots,m_n$ respectively.
Then
\begin{equation}
\label{eq:centmulti}
\bC(\cL) = k[(z_1^{m_1})^{\pm 1},\dots,(z_n^{m_n})^{\pm 1}],
\end{equation}
and so $C(\cL)$ is isomorphic to the algebra of Laurent
polynomials in $n$-variables over~$k$.
\end{corollary}

\begin{proof}  Recall (using the notation of Example \ref{ex:multiloop}) that
\[\cL = \bigoplus_{(i_1,\dots,i_n)\in \bbZ^n}
\cA_{\modi_1,\dots,\modi_n}\otk z_1^{i_1}\dots z_n^{i_n},\]
and so the inclusion ``$\supset$'' in \eqref{eq:centmulti} is clear.
For the inclusion ``$\subset$'', let $u\in \bC(\cL)$.
Now $\Stn$ is naturally $\bbZ^n$-graded and it is clear
that $\bC(\cL)$ is a graded subalgebra.  Hence we can
assume that $u = z_1^{j_1}\dots z_n^{j_n}$, where
$(j_1,\dots,j_n)\in \bbZ^n$.
But then $\cA_{\modi_1,\dots,\modi_n} \subset
\cA_{\modi_1+\modj_1,\dots,\modi_n+\modj_n}$ for all
$(i_1,\dots,i_n)\in \bbZ^n$ and   so $(\modj_1,\dots,\modj_n) = (\bar 0,\dots,\bar 0)$.
\end{proof}

\section{Untwisting iterated loop algebras}
\label{sec:untwist}

In this section we show that any $n$-step loop algebra based on a
pfgc algebra can be untwisted by an extension of the centroid of
$\cL$ that is free of finite rank.

Suppose again throughout the section that \textit{$\cL =
L(\cA,\Sigma_1,\dots,\Sigma_n)$ is an $n$-step loop algebra
based on an algebra $\cA$ over $k$}.  We use the notation of the
previous section.

It will be convenient to work with $\bC(\cL)$ rather than $C(\cL)$
(although one could use $\gamma_\cL$ to identify these algebras
using Theorem \ref{thm:centloopn} and avoid this distinction).
Note that $\bC(\cL)$ is a subalgebra of $\centA\otk \Stn$, and so
$\centA\otk \Stn/\bC(\cL)$ is a ring  extension.  This is the
extension that we use to untwist $\cL$.

We define
\[\omega_\cL :
\cL \ot_{\bC(\cL)} (\centA\otk\Stn) \to \cA\otk \Stn\] by
\[\omega_\cL(x\ot u) = u\cdot x\]
for $x\in \cL$ and $u\in \centA\otk\Stn$. One easily checks that
$\omega_\cL$ is a well-defined $\centA\otk\Stn$-algebra
homomorphism which we call canonical.

Our untwisting theorem is the following:

\begin{theorem}
\label{thm:untwist} Suppose that $\cL =
L(\cA,\Sigma_1,\dots,\Sigma_n)$ is an $n$-step loop algebra
based on a pfgc algebra $\cA$, where  $\Sigma_p$ has modulus
$m_p$ for $1\le p \le n$. Then
\begin{itemize}
\item[(i)]
$\centA\otk\Stn$ is a free $\bC(\cL)$-module of rank $m_1\dots
m_n$ with basis
\[\set{\one_{\centA} \ot z^\bi \suchthat \bi\in I_n}.\]
\item[(ii)]
The canonical map $\omega_\cL$ is an isomorphism and so
\begin{equation}
\label{eq:untwist} \cL \ot_{\bC(\cL)} (\centA\otk\Stn)
\simeq_{\centA\otk\Stn} \cA\otk \Stn.
\end{equation}
\end{itemize}
\end{theorem}

\begin{proof}
(i) Since $\bC(\cL)$ is an $n$-step loop algebra based on $\centA$
by Theorem \ref{thm:centloopn}, statement (i) follows from
Corollary \ref{cor:freemod}.

(ii) First $\cL$ spans $\cA\otk\Stn$ over $\Stn$ by Lemma
\ref{lem:canform}, and so $\cL$ spans $\cA\otk \Stn$ over
$\centA\ot \Stn$.  Thus $\omega_\cL$ is surjective.

To show that $\omega_\cL$ is injective, let $x\in
\ker(\omega_\cL)$.  Then, in particular, $x$ is an element of
$\cL \ot_{\bC(\cL)} (\centA\otk\Stn)$. Now since $\bC(\cL)$ is an
$n$-step loop algebra based on $\centA$, it follows from Lemma
\ref{lem:canform} that every element of $\centA\otk\Stn$ can be
written as a sum of elements of the form $z^\bi\cdot u$, where
$\bi\in I_n$ and $u \in \bC(\cL)$. But $z^\bi\cdot u = u
\cdot(\one_{\centA}\ot z^\bi)$. Thus from the balanced property in
the tensor product $\cL \ot_{\bC(\cL)} (\centA\otk\Stn)$, it
follows that $x$ can be written in the form
\[x= \sum_{\bi\in I_n} x_\bi \ot(\one_{\centA}\ot z^\bi),\]
where $x_\bi\in \cL$ for all $\bi$. Applying $\omega_\cL$ to this
expression yields $\sum_{\bi\in I_n} z^\bi\cdot x_\bi = 0$, and so
$x_\bi = 0$ for all $\bi\in I_n$ by Lemma \ref{lem:canform}.  Thus
$x=0$ and $\omega_\cL$ is injective.
\end{proof}

\begin{remark}
\label{rem:untwist} Suppose that $\cL$ is an $n$-step loop algebra
based on a pfgc algebra~$\cA$.

(i) We can use the canonical isomorphism $\gamma_\cL : \bC(\cL)
\to \centL$ of Theorem \ref{thm:centloopn} to identify the
algebras $\bC(\cL)$ and $\centL$. This identification is
compatible with the actions of these algebras on $\cL$ and it
gives $\centA\otk \Stn$ the structure of a $\centL$-module.  Then
\eqref{eq:untwist} can be restated as:
\begin{equation*}
\cL \ot_{\centL} \big(\centA\otk \Stn\big) \simeq_{\centA\otk
\Stn} \cA \otk \Stn. \tag{\ref{eq:untwist}$'$}
\end{equation*}
Since $\cA \otk \Stn$ is the untwisted $n$-step loop algebra based
on $\cA$, Theorem \ref{thm:untwist}  tells us that \textit{$\cL$
is untwisted by a free base ring extension of rank $m_1\dots m_n$
of the centroid of $\cL$}.

(ii) Also observe that the algebras $\cA \otk \Stn$  and $\cA
\ot_{\centA}\big(\centA\otk \Stn\big)$ are canonically isomorphic
as $\centA\otk \Stn$-algebras.  Thus the isomorphism
(\ref{eq:untwist}$'$) can be further restated as
\begin{equation*}
\cL \ot_{\centL} \big(\centA\otk \Stn\big) \simeq_{\centA\otk
\Stn} \cA \ot_{\centA}\big(\centA\otk \Stn\big).
\tag{\ref{eq:untwist}$''$}
\end{equation*}
\end{remark}

\bigskip

Theorem \ref{thm:untwist}
can  be used to
compare properties of an
iterated loop algebra as a module or algebra over its centroid
with corresponding properties of the base algebra over its
centroid. We now indicate an example of this sort
of argument.

\begin{corollary}
\label{cor:CL} Let $\cL$  be an $n$-step loop algebra based
on a pfgc algebra $\cA$. If $\cA$ is a projective $\centA$-module
then $\cL$ is a finitely generated projective $C(\cL)$-module
\end{corollary}

\begin{proof} As in Remark \ref{rem:untwist}, we identify
$\bC(\cL)$ and $\centL$ using $\gamma_\cL$. By axiom P2 of pfgc
algebras and the present assumption, $\cA$ is a finitely generated
projective $\centA$-module. Hence $\cA\ot_{\centA}
\big(\centA\ot_k S^{\ot n}\big)$ is a finitely generated
projective $\centA\ot_k S^{\ot n}$-module. Thus by
(\ref{eq:untwist}$''$), $\cL \ot_{\centL} \big(\centA\otk
\Stn\big)$ is a finitely generated projective $\centA\ot_k S^{\ot
n}$-module.  But the extension $\centA\ot_k S^{\ot n}/C(\cL)$ is
free of finite rank by Theorem \ref{thm:untwist}(i), and so it is
faithfully flat. The result  now follows from \cite[Ch.~I, \S~3.6,
Prop.~12]{B:CA}
\end{proof}

In  the same spirit, we now describe an application
of Theorem  \ref{thm:untwist} for associative algebras. For this purpose we first
recall some definitions and basic facts about Azumaya algebras. A unital
associative algebra $\cD$ over a
ring $B$ is called an \emph{Azumaya algebra} over $B$ if $\cD$ is
central and separable over $B$ (see for example \cite[\S 5]{KO}).
If $\cD$ is an Azumaya algebra over $B$, then $\cD$ is a finitely
generated projective $B$-module \cite[Th\'eor\`eme~5.1]{KO}, and
so $\cD_{\mathfrak m}$ is a free $B_{\mathfrak m}$-module of
finite rank $r_\mathfrak m$ for each maximal ideal $\mathfrak m$
of $B$.  $\cD$ is said to have \textit{constant rank} $r$ over $B$
if $r_\mathfrak m = r$ for all such $\mathfrak m$ \cite[\S
II.5.3]{B:CA}.  It is known that if $\cD$
is a unital associative algebra over a ring $B$ and $\ell$ is a positive integer then
\begin{equation}
\label{eq:Azumaya}
\parbox{4in}
    {$\cD$ is an Azumaya algebra of constant rank $\ell^2$ over $B$
    if and only if
    there exists a faithfully flat extension $B'/B$ of
    rings so that $\cD\ot_B B'$ is isomorphic as a $B'$-algebra to
    the algebra $M_\ell(B')$ of $\ell\times\ell$-matrices
    over $B'$.}
\end{equation}
In that case we will say that
$\cD$ is \textit{split} by the extension $B'/B$.
Indeed the implication ``$\Rightarrow$'' of
\eqref{eq:Azumaya} is Corollary 6.7 of \cite{KO}.
For the converse, the algebra $M_\ell(B')$ is an Azumaya algebra of constant
rank $\ell^2$ over $B'$, and hence $\cD$ is an Azumaya algebra of constant rank $\ell^2$
over $B$ since the extension $B'/B$ is faithfully flat (see Lemma 5.1.9(1) in \cite{Kn}
and Exercise 8 in \cite[\S~II.5]{B:CA}).

\begin{corollary}
\label{cor:Azumaya}
Suppose that
$\cL$ is an $n$-step loop algebra
based on the matrix algebra $M_\ell(k)$ over $k$.
Then $\cL$ is a prime Azumaya algebra of finite rank $\ell^2$
over its centroid $C(\cL)$ which is split by the extension
$\Stn/C(\cL)$.
\end{corollary}

\begin{proof}
Let $\cA = M_\ell(k)$.  Then $\cA$ is a  prime unital associative
algebra over $k$ and hence so is $\cL$ (by Theorem \ref{thm:basic}).
Also $C(\cA) = k$ and so $C(\cA)\otk \Stn = \Stn$ as in
Remark \ref{rem:centcalc}.  Thus by \eqref{eq:untwist}$'$ we have
\[
\cL \ot_{\centL} \Stn \simeq_{\Stn} M_\ell(k) \otk \Stn
\simeq_{\Stn} M_{\ell}(\Stn).\] Our conclusion now follows from
\eqref{eq:Azumaya}, since the  extension $\Stn/C(\cL)$ is
faithfully  flat.
\end{proof}

\section{Permanence of  type}
\label{sec:permanence}

There is a classical notion of   type for
simple pfgc Lie algebras in characteristic zero (see Example
\ref{ex:typeLie} below). This notion can easily be extended using
the central closure to include prime pfgc Lie algebras in
characteristic $0$. It will be a consequence of the results in
this section that  type is preserved under
the loop construction (that is  type is
permanent).

An analogous notion of  type can be
defined for many other important classes of prime pfgc algebras
besides Lie algebras in characteristic 0.  Moreover, since
algebras in these classes arise naturally as coordinate algebras
in the study of Lie algebras, and in particular in the study of
extended affine Lie algebras, it is desirable to include these
classes in our discussion of  type. This
generality requires almost no extra effort once the appropriate
definitions are made.  That being said, the reader can safely
choose to assume throughout this section that the base algebras,
and hence their loop algebras, are Lie algebras in characteristic
$0$.

We begin by recalling the classical notion of
type for simple pfgc Lie algebras in characteristic $0$.

\begin{example}
\label{ex:typeLie} Suppose that $k$ has characteristic $0$.  Let
$\cA$ be a simple pfgc Lie algebra over $k$ (or equivalently let
$\cA$ be a simple Lie algebra over $k$ that is finitely generated
as a module over its centroid). Then, since $\cA$ is simple, it is
easily checked that $\centA$ is a field.  Hence, if we let $K$ be
an algebraic closure of $\centA$, the algebra  $\cA \ot_{\centA}
K$  is a finite dimensional simple Lie algebra over $K$ by Lemma
\ref{lem:centalsimple}.  The \textit{type} of $\cA$ is defined in
\cite[Ch.~X, \S~3]{J1} to coincide with the type of the root
system of the $K$-algebra $\cA \ot_{\centA} K$ relative  to any
Cartan subalgebra.
\end{example}

In order to extend this notion to other classes of algebras, we need to
introduce some  terminology.

\begin{definition}
Recall that a \textit{variety} over $k$ is a class $\bV$ of
algebras over $k$ that is defined by a set of identities in the
free nonassociative algebra $k_\text{na}[x_1,x_2,\dots]$ in
countably many symbols \cite[\S~1.2]{ZSSS}. A variety $\bV$ over
$k$ is said to be \textit{homogenous} if the ideal in
$k_\text{na}[x_1,x_2,\dots]$ consisting of all identities
satisfied by all algebras in $\bV$ is homogeneous. Algebras in a
variety $\bV$ over $k$ will be simply called
\textit{$\bV$-algebras}.
\end{definition}

A very familiar example is the variety $\bV$ of Lie algebras over
$k$ which is defined by the identities  $x_1x_1$ and $(x_1x_2)x_3 +
(x_2x_3)x_1 + (x_3x_1)x_2$. In that case $\bV$  is homogeneous
\cite[\S~1.4]{ZSSS}, and a $\bV$-algebra is just a Lie algebra
over $k$.

\begin{remark}  Suppose that
$\bV$ is a \textit{variety} over $k$.  Suppose that $B$ is a
unital associative commutative $k$-algebra. A
\textit{$\bV$-algebra over $B$} will mean a $B$-algebra $\cA$ with
the property that $\cA$, when regarded as an algebra over $k$, is
in $\bV$.  In other words, a $\bV$-algebra over $B$ is a
$B$-algebra that satisfies the identities (which are identities
with coefficients from our fixed base field $k$) that define
$\bV$.
\end{remark}

The  homogeneity  assumption is important for our purposes since
homogeneous varieties are closed under base ring extension.

\begin{lemma}
\label{lem:varietytensor} Suppose that $\bV$ is a homogeneous
variety over $k$ and $B\to B'$ is a homomorphism of unital
commutative associative $k$-algebras.  If $\cA$ is a $\bV$-algebra
over $B$ then $\cA \ot_B B'$ is a $\bV$-algebra over $B'$.
\end{lemma}
\begin{proof}
This follows the proof of Theorem~6 in \cite[\S~1.4]{ZSSS}.
\end{proof}

\begin{corollary}
\label{cor:varietytensor} Suppose that $\bV$ is a homogeneous
variety over $k$.  If $\cL$ is an $n$-step loop algebra based on a
$\bV$-algebra $\cA$, then $\cL$ is a $\bV$-algebra.
\end{corollary}

\begin{proof} By Lemma \ref{lem:varietytensor}, $\cA \otk \Stn$ is a
$\bV$-algebra. Hence so is its  subalgebra~$\cL$.
\end{proof}

We will be interested in homogeneous varieties $\bV$ that satisfy
the following axiom:
\begin{itemize}
\item[(S)] If $K/k$ is a field extension and $\cA$ is a finite
dimensional semiprime $\bV$-algebra over $K$ then $\cA$ is a
direct sum of simple algebras over $K$.
\end{itemize}

\begin{example}
\label{ex:axiomS} In each of the following cases, the variety
$\bV$ is homogeneous and satisfies axiom (S):
\begin{itemize}
    \item[(a)] $\charr(k)=0$, $\bV$ is the variety of Lie algebras.
    \item[(b)] $\bV$ is the variety of associative algebras.
    \item[(c)] $\bV$ is the variety of commutative associative algebras.
    \item[(d)] $\bV$ is the variety of alternative algebras.
    \item[(e)] $\charr(k)\ne 2$, $\bV$ is the variety of Jordan algebras.
\end{itemize}
Indeed the fact that these varieties are homogeneous is proved in
\cite[\S~1.4]{ZSSS}.  Moreover axiom (S) follows from the
structure theory for the variety $\bV$ in each case. For  example,
in case (a), suppose that $K/k$ is a field extension and $\cA$ is
a finite dimensional semiprime Lie algebra over $K$. If the
radical $\mathcal R$ of $\cA$ is nonzero, then the last nonzero
term in the derived series for $\mathcal R$ has trivial
multiplication, contrary to the assumption that $\cA$ is
semiprime.  So $\mathcal R = 0$ and hence $\cA$ is the direct sum
of simple algebras \cite[\S~III.4]{J1}. Similarly we can use (for
example) \cite[\S 12.2, Theorem~3]{ZSSS} in cases (b), (c) and (d)
and \cite[\S~V.2, Lemma~2 and \S~V.5, Corollary~2]{J2} in case (e)
to verify axiom (S).
\end{example}

The reason for our interest in Axiom (S) is that it allows us to
prove the following proposition.

\begin{proposition}
\label{prop:cquot2} Let $\bV$ be a homogeneous variety over $k$
that satisfies axiom {\rm(S)}.  Suppose that $\cA$ is a prime pfgc
$\bV$-algebra over $k$ and let $\tcentA$ be the quotient field of
$\centA$. Then the central closure $\tA := \cA\ot_{\centA}
\tcentA$ of $\cA$ is a finite dimensional central simple
$\bV$-algebra over~$\tcentA$.
\end{proposition}

\begin{proof}  By Proposition~\ref{prop:cquot1}, $\tA$ is a prime pfgc algebra
over $k$ and hence also over $\tcentA$ (see Remark
\ref{rem:pfgcinvariant} and Lemma \ref{lem:primeinvariant}). Also,
by Lemma \ref{lem:varietytensor}, $\tA$ is a $\bV$-algebra. Hence,
by axiom (S), $\tA$ is the direct sum of simple algebras over
$\tcentA$.  Since $\tA$ is prime, there is only one summand in
this sum. Thus $\tA$ is  a simple algebra over $\tcentA$.  Finally, by
Proposition~\ref{prop:cquot1}, $\tA$ is central and finite  dimensional
over $\tcentA$.
\end{proof}

\begin{remark}
\label{rem:affine}
Suppose that $\cA$ is as in Proposition \ref{prop:cquot2}.
Then in the  terminology of \cite[\S~1]{PS}, Proposition \ref{prop:cquot2}
says that $\cA$ is a \emph{central order} in the finite dimensional
central simple algebra $\tA$.
\end{remark}

We will also need a set $\bX$ of algebras over $k$ that play the
role of the split simple Lie algebras over $k$.

\begin{definition}
\label{def:archetype} Suppose that $\bV$ is a homogenous variety
over $k$.  A \textit{set of archetypes} for $\bV$ is a set $\bX$
of finite dimensional central simple $\bV$-algebras over $k$ such
that the following axioms hold:
\begin{itemize}
\item[(A1)]
If $K/k$ is an algebraically closed field extension and $\cA$ is a
finite dimensional central simple $\bV$-algebra over $K$ then
there exists  $\cX\in \bX$ so that $\cA \simeq_K \cX\otk K$.
\item[(A2)] If $K/k$ is a field extension and $\cX,\cX'\in \bX$
then
\[\cX\otk K\simeq_K\cX'\otk K \implies
\cX = \cX'\]
\end{itemize}
In particular, the elements of $\bX$ are pairwise nonisomorphic
over $k$.
\end{definition}

\begin{example}
\label{ex:archetype} In each of the cases (a)-(e) in Example
\ref{ex:axiomS} there is a natural choice for a set $\bX$ of
archetypes:
\begin{itemize}
    \item[(a)] $\charr(k)=0$, $\bV$ is the variety of Lie algebras and
    $\bX = \set{\cX_\Pi}$, where $\Pi$ runs through all connected
    Dynkin diagrams (up to isomorphism) and $\cX_\Pi$ denotes the split
    simple Lie algebra with Dynkin diagram $\Pi$
    \cite[\S~VII.4]{J1}.
    \item[(b)] $\bV$ is the variety of associative algebras and
    $\bX = \set{\Mat_\ell(k) \suchthat \ell\ge 1}$, where
    $\Mat_\ell(k)$ is the algebra of $\ell\times\ell$-matrices over $k$.
    \item[(c)] $\bV$ is the variety of commutative associative algebras and
    $\bX = \set{k}$.
    \item[(d)] $\bV$ is the variety of alternative algebras
    and $\bX = \set{\Mat_\ell(k) \suchthat \ell\ge 1}\cup\set{\mathcal O}$,
    where $\mathcal O$ denotes the split Cayley-Dickson (= octonion) algebra
    \cite[\S~2.4]{ZSSS}.
    \item[(e)] $\charr(k)\ne 2$, $\bV$ is the variety of Jordan algebras
    and $\bX$ is the set consisting of the following algebras: $k$;
    the Jordan algebra  constructed from a
    nondegenerate symmetric bilinear form
    with matrix $\diag(1,-1,\dots,1,-1)$ on a $2\ell$-dimensional space
    over $k$, $\ell \ge 1$;
    the Jordan algebra constructed from
    a nondegenerate symmetric bilinear form
    with matrix $\diag(1,-1,\dots,1,-1,1)$ on a $2\ell+1$-dimensional space over $k$,
    $\ell \ge 1$;
    the algebra of $\ell\times \ell$ hermitian matrices with
    coordinates from the split composition
    algebras of dimension 1, 2 and 4, $\ell\ge 3$;
    and the algebra of $3\times 3$ hermitian matrices
    over $\mathcal O$ \cite[\S~1.4 and~4.3]{J2}.
    \end{itemize}
The fact that $\bX$ satisfies axioms (A1) and (A2) follows from
the classification of finite dimensional central simple algebras
over algebraically closed fields in each case.  See for example
\cite[\S~IV.3, Theorem~3]{J1} in case (a), \cite[\S~12.2,
Theorem~3]{ZSSS} in cases (b), (c) and (d), and \cite[\S~V.6,
Corollary~2]{J2} in case (e).
\end{example}

\begin{remark}
\label{rem:alternate} A homogeneous variety $\bV$ over $k$ may
possess more than one set of archetypes.  For example if $k =
\mathbb{R}$ and $\bV$ is the variety of Lie algebras over $k$ as
in Example \ref{ex:archetype}(a) above, then an alternate choice
of a set of archetypes is the set $\bX = \set{\mathcal C_\Pi}$,
where $\Pi$ runs through all connected Dynkin diagrams (up to
isomorphism) and $\mathcal C_\Pi$ denotes the compact real Lie
algebra whose complexification is the simple Lie algebra with
Dynkin diagram $\Pi$.
\end{remark}

\noindent \textbf{Assumption:} \textit{For  the rest of this
section we assume that $\bV$ is a homogeneous variety over $k$
that satisfies axiom {\rm(S)}, and  that there exists a set $\bX$
(which we fix) of archetypes for $\bV$.}

\bigskip

We can now prove the proposition that allows us to define the
notion of  type.

\begin{proposition}
\label{prop:typeinv} Suppose that $\cA$ is a prime pfgc
$\bV$-algebra over $k$. If
\[\centA \hookrightarrow K\]
is a unital $k$-algebra monomorphism of $\centA$ into an
algebraically closed field extension $K$ of $k$  (such a
monomorphism exists since $C(\cA)$ is an integral domain),  then
there exists a unique $\cX\in \bX$ so that
\begin{equation}
\label{eq:type}
    \cA \ot_{\centA} K \simeq_K \cX \otk K,
\end{equation}
where on the left $K$ is regarded as an algebra over $\centA$
using the given monomorphism. Moreover, $\cX$ is independent of
the choice of the $k$-algebra monomorphism $\centA \hookrightarrow
K$.
\end{proposition}

\begin{proof}
First let $L$ be an algebraic closure of $\tcentA$. By
Proposition~\ref{prop:cquot2}, $\tA$ is a finite dimensional
central simple $\bV$-algebra over $\tcentA$. Therefore, by
Lemma~\ref{lem:centalsimple}, $\tA  \ot_{\tcentA} L$ is a finite
dimensional central simple algebra over $L$.  So, by
Lemma~\ref{lem:varietytensor}, $\tA  \ot_{\tcentA} L$ is a finite
dimensional central simple $\bV$-algebra over~$L$. Thus, by axiom
(A1) (see Definition \ref{def:archetype}), there exists $\cX\in
\bX$ so that $\tA  \ot_{\tcentA} L \simeq_L \cX\otk L$. Then
\begin{equation}
\label{eq:typeL} \cA \ot_{\centA} L \simeq_L (\cA \ot_{\centA}
\tcentA)\ot_{\tcentA} L = \tA  \ot_{\tcentA} L \simeq_L \cX\otk L.
\end{equation}

Now let $\centA \hookrightarrow K$ be an arbitrary unital
$k$-algebra monomorphism of $\centA$ into an algebraically closed
extension $K$ of $k$.  This extends to a unital $k$-algebra
monomorphism $\tcentA \hookrightarrow K$ which in turns extends to
a unital $k$-algebra monomorphism $L \hookrightarrow K$.  We use
this latter monomorphism to identify $L$ as a subfield of $K$.
Then using \eqref{eq:typeL} we have
\[ \cA \ot_{\centA} K \simeq _K (\cA \ot_{\centA} L)\ot_L K
\simeq_K (\cX\otk L)\ot_L K \simeq_K \cX \otk K.
\]
This shows the existence of an element $\cX\in \bX$ satisfying
\eqref{eq:type}.  The uniqueness follows from Axiom (A2).

Finally if $\centA \hookrightarrow K'$ is another unital
$k$-algebra monomorphism of $\centA$ into an algebraically closed
extension $K'$ of $k$,  then the argument just given using
\eqref{eq:typeL} shows that $\cA \ot_{\centA} K' \simeq_{K'} \cX
\otk K'$.
\end{proof}

\begin{definition}
\label{def:type} Let $\cA$ be a prime pfgc $\bV$-algebra over $k$.
The element $\cX\in \bX$ described in
Proposition~\ref{prop:typeinv} is called the \textit{type}
of $\cA$ (relative to $\bX$) and denoted by
$t(\cA)$.  We also
sometimes refer to $t(\cA)$ as the \textit{absolute type}
of $\cA$ since it is determined by extending the base ring $C(\cA)$ to an algebraically
closed field.

\end{definition}

\begin{example}
\label{ex:typeLie2} Let $\charr(k)=0$, let $\bV$  be the variety
of Lie algebras, and let $\bX = \set{\cX_\Pi}$ be as in Example
\ref{ex:archetype}(a). If we identify $\cX_\Pi$ with the diagram
$\Pi$, then Definition \ref{def:type} assigns to each prime pfgc
$\bV$-algebra $\cA$ a connected Dynkin diagram $t(\cA)$. (If $\cA$
is a simple pfgc algebra, this is exactly what was done in
 Example \ref{ex:typeLie}.)
\end{example}

The following result tells us that   type
is an isomorphism invariant for prime pfgc algebras.
\begin{proposition}
\label{prop:typeisoinv} Suppose that $\cA$ and $\cA'$ are prime
pfgc $\bV$-algebras over $k$. Then
\[ \cA \simeq_k \cA' \implies t(\cA) = t(\cA').\]
\end{proposition}

\begin{proof}
Let $\varphi : C(\cA') \hookrightarrow K$ be a unital $k$-algebra
monomorphism of $C(\cA')$ into an algebraically closed field
extension $K$ of $k$. Denote the resulting  action of $C(\cA')$ on
$K$ by $(\chi',\alpha)\mapsto \chi'\cdot\alpha$.

Fix a $k$-algebra isomorphism $\rho  :\cA\to \cA'.$ Then $\rho$
induces an isomorphism $C(\rho) : \centA\to C(\cA')$ by Lemma
\ref{lem:centinduce}. So the composite map $\varphi \circ C(\rho)
: \centA \to K$ is a unital $k$-algebra monomorphism which we use
to view $K$ as an algebra over $\centA$. The resulting action of
$\centA$ on $K$ is given by
\[
\chi \cdot \alpha  = C(\rho)(\chi)\cdot \alpha.
\]
for $\chi\in \centA$ and $\alpha\in K$.

The bi-additive map $\wt \rho :\cA\times K\to \cA'\ot_{C(\cA')} K$
satisfying $ \wt\rho :(a,\alpha  )\mapsto \rho(a)\ot \alpha $ is
then $\centA$-balanced. Indeed if $\chi  \in \centA$, $a\in \cA$
and $\alpha\in K$ we have
$$
\begin{aligned}
\wt\rho \big(\chi(a),\alpha  \big) &=\rho\big(\chi(a)\big) \ot
\alpha  = C(\rho)(\chi)\big(\rho
(a)\big)\ot \alpha  \\
&= \rho(a)\ot \big(C(\rho)(\chi)\cdot \alpha\big) =\rho(a)\ot \chi
\cdot \alpha = \wt \rho (a,\chi  \cdot \alpha).
\end{aligned}
$$
Thus $\wt\rho$ induces a $k$-linear map $\cA\ot_{\centA} K \to
\cA'\ot_{C(\cA')} K$ so that $a\ot\alpha \mapsto \rho(a) \ot
\alpha$ for $a\in \cA$ and $\alpha\in K$.  This map is clearly a
homomorphism of $K$-algebras. In a similar fashion we obtain a
homomorphism of $K$-algebras $\cA'\ot_{C(\cA')} K\to
\cA\ot_{\centA} K$ so that $a' \ot \alpha \mapsto
\rho^{-1}(a')\ot\alpha$ for $a'\in \cA'$  and $\alpha\in K$.
These maps are inverses of each other and so we  have
\[\cA\ot_{\centA}K\simeq_{K} \cA'\ot_{C(\cA')} K.\]
Thus, $\cX\otk K\simeq_{K} \cX'\otk K$, where $\cX = t(\cA)$ and
$\cX' = t(\cA')$, and so $t(\cA)=t(\cA')$.
\end{proof}

Our main result in this section is the following:

\begin{theorem}
\label{thm:permanence} \emph{(Permanence of type)} If   $\cL$ is an $n$-step loop algebra based on a
prime pfgc $\bV$-algebra $\cA$, then $\cL$ is a prime pfgc
$\bV$-algebra and
\[t(\cL) = t(\cA).\]
\end{theorem}

\begin{proof}
By Theorem~\ref{thm:basic}(iii) and  (iv) and Corollary
\ref{cor:varietytensor}, $\cL$ is a prime pfgc $\bV$-algebra. So
$t(\cA)$ and $t(\cL)$ are defined and it remains so show that
these   types  are equal. For this note
first that $\centA$ is an integral domain by
Lemma~\ref{lem:intdom}(i), and so the algebra $\centA\otk \Stn
\simeq_k \centA[z_1^{\pm 1},\dots,z_n^{\pm 1}]$ is an integral
domain. Let $K$ be an algebraic closure of the quotient field  of
$\centA\otk \Stn$. Now by (\ref{eq:untwist}$''$) we have the
isomorphism
\[\cL \ot_{\centL} \big(\centA\otk \Stn\big)
\simeq_{\centA\otk \Stn} \cA \ot_{\centA}\big(\centA\otk
\Stn\big).\] Tensoring this over $\centA\otk \Stn$ with $K$ yields
the isomorphism
\[\cL \ot_{\centL} K  \simeq_{K}
\cA \ot_{\centA} K.\] Hence we have $\cX \otk K \simeq_K \cX' \otk
K$, where $\cX = t(\cA)$ and $\cX' = t(\cL)$, and so $\cX = \cX'$.
\end{proof}

Since finite dimensional central simple algebras over $k$ are
prime pfgc algebras, we have:

\begin{corollary}
\label{cor:centsimplebaseb} If $\cL$ is an $n$-step loop algebra
based on a finite dimensional central simple $\bV$-algebra $\cA$ over $k$,
then $\cL$ is a prime pfgc $\bV$-algebra and $t(\cL) = t(\cA)$.
\end{corollary}

\begin{corollary}
\label{cor:archetypebase} If $\cX\in \bX$ and $\cL$ is an $n$-step
loop algebra based on $\cX$, then $\cL$ is  a prime pfgc
$\bV$-algebra of   type $\cX$. If further
$\cX'\in \bX$ and $\cL'$ is an $n'$-step loop algebra based on
$\cX'$, then
\[
\cL \simeq_k \cL' \Longrightarrow \cX = \cX'\text{ and } n=n'.
\]
\end{corollary}

\begin{proof} The first statement follows from Corollary \ref{cor:centsimplebaseb}
since $t(\cX) = \cX$.  For the second statement suppose that $\cL
\simeq_k \cL'$. Then by Proposition \ref{prop:typeisoinv}, $t(\cL)
= t(\cL')$ and so $\cX = \cX'$.  Finally by Corollary
\ref{cor:centsimplebasea}, $n=n'$.
\end{proof}

Our primary focus in future work will be on the case when the base
algebra is a finite dimensional split simple Lie algebra. For ease
of reference we therefore record Corollary \ref{cor:archetypebase}
explicitly in that case.

\begin{corollary}
\label{cor:Liebase} Suppose that $\cA$ is a finite dimensional
split simple Lie algebra over a field $k$ of characteristic $0$,
and $\cL$ is an $n$-step loop algebra based on $\cA$.  Then $\cL$
is a prime pfgc Lie algebra and for any unital $k$-algebra
monomorphism $\centL \hookrightarrow K$ into an algebraically
closed extension $K$ of $k$ we have
\[\cL \ot_{\centL} K \simeq_K \cA \otk K.\]
Moreover, if $\cA'$  is a finite dimensional split simple Lie
algebra and $\cL'$ is an $n'$-step loop algebra based on $\cA'$,
then
\begin{equation}
\label{eq:isominvariant}
\cL \simeq_k \cL' \Longrightarrow \cA \simeq_k \cA'\text{ and } n=n'.
\end{equation}
\end{corollary}

\begin{proof}  We apply Corollary \ref{cor:archetypebase} to the
case when $\bV$ is the variety of Lie algebras and $\bX =
\set{\cX_\Pi}$ as in Example \ref{ex:archetype}(a). Since any
finite dimensional split simple Lie algebra over $k$ is isomorphic
to exactly one algebra in $\bX$ the result follows.
\end{proof}

If $\cL$  is an $n$-step loop algebra
based on a finite dimensional split simple Lie algebra $\cA$ (in
characteristic 0), then \eqref{eq:isominvariant} tells us that
both (the isomorphism class of) the base algebra $\cA$ and the
number of steps $n$ are isomorphism invariants of~$\cL$.
This  answers a natural question that began
the research described in this paper.  We have now seen in
Corollary \ref{cor:archetypebase} that the result is true in a
much broader context.  The interested reader can easily write down
the results corresponding to Corollary \ref{cor:Liebase}  for the
varieties of associative algebras, alternative algebras and Jordan
algebras (see Example \ref{ex:archetype}(b), (d) and
(e)).

\section{Two-step loop algebras}
\label{sec:twostep}

In  this section we look more closely at
2-step iterated loop algebras and their centroids. We then conclude with
a detailed look at two examples that illustrate several of the concepts studied in
this article.

\textit{Throughout this section, we assume that $m_1$ and $m_2$ are positive integers
and that $k$ contains a primitive $m_i^\thsup$-root of unity $\zeta_{m_i}$, $i=1,2$.}
We use the notation of Section \ref{sec:iterated} (for iterated loop algebras).

We start with some further notation.  Let $k^\times = \set{\rho\in k \suchthat \rho \ne 0}$
be the group of units of $k$.  If $\rho\in k^\times$, we   let
\[\strange\]
denote the unital associative commutative $k$-algebra presented by the generators
$u_1$, $u_2$, $u_2^{-1}$, $w$ subject to the relations
\[u_2 u_2^{-1} = 1 \andd w^2 = (u_1^2 - 4\rho )u_2.\]
It is clear that the set
\[\set{u_1^{i_1} u_2^{i_1} w^j :  i_1 \in \bbZ_{\ge 0},\  i_2 \in \bbZ,\ j = 0,1}\]
is a $k$-basis for $\strange$. It is  also easy to verify that the
group of units of $\strange$ is given by
\begin{equation}
\label{eq:strangeunit}
U(\strange) = \set{\alpha u_2^{i_2} \suchthat \alpha\in k^\times,\ i_2\in \bbZ}.
\end{equation}
Indeed, one way to see this is to make use of the multiplicative norm function
$N: \strange \to k[u_1,u_2^{\pm 1}]$ defined by $N(a_1 + a_2 w) = a_1^2 - a_2^2 w^2$ for
$a_1,a_2 \in k[u_1,u_2^{\pm 1}]$, and use the fact that if $u$ is a unit
in $\strange$ then $N(u)$ is a unit in $k[u_1,u_2^{\pm 1}]$.  We leave the details
of this to the reader.

The algebra $\strange$ is important in the study of iterated loop algebras because
of the following fact.

\begin{lemma}
\label{lem:kind}
Let
$\cL = L(k,\Sigma_1,\Sigma_2)$ be a $2$-step iterated loop algebra based on the algebra
$k$, where $\Sigma_i$ has modulus $m_i$ for $i=1,2$.  Then exactly one of the following holds:
\begin{itemize}
    \item[(a)] $\cL \simeq_k k[t_1^{\pm 1},t_2^{\pm 2}]$
    (the algebra of Laurent polynomials in 2 variables).
    \item[(b)] $\cL \simeq_k \strange$ for some $\rho\in k^\times$.
\end{itemize}
Moreover (a) holds if and only if $z_1^{m_1} z_2^j \in \cL$ for some $j\in \bbZ$.
\end{lemma}

\begin{proof}
Note that the group of units in $k[t_1^{\pm 1},t_2^{\pm 1}]$ spans
the algebra $k[t_1^{\pm 1},t_2^{\pm 1}]$, whereas this is not true
for the algebra $\strange$ (by \eqref{eq:strangeunit}). Thus (a)
and (b) cannot hold simultaneously. So it  remains to show that
either (a) or (b) holds (the final statement will be proved along
the way).

Now as noted in Remark \ref{rem:iteratedloop}(iii), we have
\[ \cL = L(k,\sigma_1,\sigma_2),\]
where $\sigma_1$ is an automorphism of period $m_1$ of $\cL_0 = k$,
and $\sigma_2$ is an automorphism of period $m_2$ of $\cL_1 = L(k,\sigma_1)$.
Then, since $\sigma_1$ is an algebra homomorphism, $\sigma_1 = 1$ and so
\[\cL_1 = k[y_1^{\pm 1}], \text{  where  } y_1 = z_1^{m_1}.\]
Thus $\sigma_2$ is an automorphism of period $m_2$ of $k[y_1^{\pm
1}]$. Hence either $\sigma_2(y_1) = \rho y_1$ for some
$m_2^\thsup$ root of unity $\rho$ in $k^\times$ or $\sigma_2(y_1)
= \rho y_1^{-1}$ for some $\rho\in k^\times$. Moreover (for the
proof of the final statement in the proposition) the first of
these possibilities holds if and only if $y_1$ is homogeneous in
the grading $\Sigma_2$ determined by $\sigma_2$ which in turn
holds if and only if $y_1z_2^j\in \cL$ for some $j\in \bbZ$.

\textit{Case (a)}:  Suppose that $\sigma_2(y_1) = \rho y_1$ for some $m_2^\thsup$ root of unity
$\rho$ in $k^\times$.  Let $n_2$ be the order of $\rho$ in $k^\times$.
Then $n_2$ is a divisor of $m_2$,
\[\rho = \zeta_{m_2}^{p_2r}, \text{  where  } p_2 = \frac {m_2}{n_2},\]
and $r$  is relatively prime to $n_2$ (take $r=0$ if $n_2 = 1$).
Choose an inverse $s$ for $r$ modulo $n_2$ (take $s=0$ if $n_2 = 1$).
Now the grading $\Sigma_2$ of $\cL_1$ is given by
$\cL_1 = \bigoplus_{\modj\in \Zmod{m_2}} (\cL_1)_{\modj}$, where
$(\cL_1)_{\modj}$ is spanned by the elements $y_1^i$ with
$i\in \bbZ$ and $\sigma_2(y_1^i) = \zeta_{m_2}^j y_1^i$.  But
$n_2$ and $s$ are relatively prime and so any integer can be expressed
in the form $an_2 + bs$, where $a,b\in \bbZ$.  Also
\[\sigma_2(y_1^{an_2+bs}) = \rho^{an_2+bs}y_1^{an_2+bs}  = \rho^{bs}y_1^{an_2+bs}
= \zeta_{m_2}^{p_2rbs}y_1^{an_2+bs} =
\zeta_{m_2}^{p_2b}y_1^{an_2+bs}\] and so $y_1^{an_2+bs}\in
(\cL_1)_{\overline{p_2b}}$.  Therefore $\cL = L(\cL_1,\sigma_2)$
is spanned by elements of the form
\[
y_1^{an_2+bs}z_2^{p_2 b},
\quad a,b\in \bbZ.
\]
 But $y_1^{an_2+bs}z_2^{p_2 b} =
(y_1^{n_2})^a(y_1^s z_2^{p_2})^b$. Hence we obtain
\[\cL = k[t_1^{\pm 1},t_2^{\pm 2}], \text{  where } t_1 = y_1^{n_2} \text{ and } t_2 = y_1^s z_2^{p_2}.\]

\textit{Case (b)}: Suppose that $\sigma_2(y_1) = \rho y_1^{-1}$ for some $\rho\in k^\times$.
Then $\sigma_2$ has order $2$ and so $m_2$ is even.  Let $p_2 = \frac {m_2} 2$
and $y_2 = z_2^{p_2}$.  Then
\[\cL = \left(\cL_1^+\ot_k k[(y_2^2)^{\pm 1}]\right)\oplus
\left(\cL_1^-\ot_k y_2k[(y_2^2)^{\pm 1}]\right), \]
where $\cL_1^{\pm}$ is the $\pm 1$-eigenspace for $\sigma_2$.
Now it is clear that $\cL_1^+$ has a $k$-basis consisting of the elements $(y_1+\rho y_1^{-1})^a$,
$a\ge 0$.   Therefore
$\cL_1^+\ot_k k[(y_2^2)^{\pm 1}]$ has basis
\[(y_1+\rho y_1^{-1})^a y_2^{2b}, \quad a,b\in \bbZ,\ a\ge 0.\]
Also  one easily checks that $\cL_1^- = (y_1-\rho y_1^{-1})\cL_1^+$, and so
$\cL_1^-\ot_k y_2k[(y_2^2)^{\pm 1}]$ has basis
\[(y_1-\rho y_1^{-1})y_2(y_1+\rho y_1^{-1})^a y_2^{2b}, \quad a,b\in \bbZ,\ a\ge 0.\]
Thus, setting
\[u_1 = y_1+\rho y_1^{-1},\ u_2 = y_2^2 \andd w = (y_1-\rho y_1^{-1})y_2,\]
we see that $\cL$ has basis
$u_1^{a} u_2^{b} w^c$,
$a \in \bbZ_{\ge 0}$, $b \in \bbZ$, $c = 0,1$.  Moreover, one checks directly that
$w^2 = (u_1^2 - 4\rho )u_2$, and so we have identified $\cL$ with $\strange$.
\end{proof}

\begin{remark} Lemma \ref{lem:kind}(a)
is an immediate consequence of a more general ``erasing theorem''
that was proved in \cite[Theorem 5.1]{ABP2}.  We have included the proof above since
it is short and self contained.
\label{rem:altapproach1}
\end{remark}

\begin{remark}
\label{rem:altapproach2}   If $\rho, \rho'\in k^\times$,
one can show that
\[\strange \simeq_k \strangep \iff \rho' \rho^{-1} \text{ is a square in } k^\times.\]
In particular, if $k$ is algebraically closed, the isomorphism
class of $\strange$ does not depend on $\rho$. In that case Lemma
\ref{lem:kind} tells us that, up to isomorphism, there are exactly
two (one step) loop algebras based on $k[y_1^{\pm 1}]$.  This fact
is a special case of a more general result about (one step) loop
algebras based on the algebra $\cA$ of Laurent polynomials
$k[y_1^{\pm 1},\dots,y_q^{\pm 1}]$ over an algebraically closed
field~$k$. Indeed, using the fact that the abstract  automorphism
group of $\cA$ is $(k^\times)^q\times \GL_q(\bbZ)$ and some
techniques from Galois cohomology (see Remark \ref{rem:untwist1}),
one can show that there is an injective map that attaches to each
$R$-isomorphism class of loop algebras based on $\cA$ an invariant
in the set of conjugacy classes of $GL_q(\bbZ)$. (When $q=1$,
$GL_q(\bbZ)$ has exactly two conjugacy classes and  one can show
that $R$-isomorphism coincides with $k$-isomorphism.) We omit
proofs of the statements in this remark, since we will not be
using these statements here and since their proofs would take us
rather far  afield.
\end{remark}

Lemma \ref{lem:kind} together with Theorem \ref{thm:centloopn}
implies the following more general result:

\begin{proposition}
\label{prop:kind}
Let
$\cL = L(\cA,\Sigma_1,\Sigma_2)$ be a $2$-step iterated loop algebra based on a finite dimensional
central simple algebra $\cA$ over $k$,
where $\Sigma_i$ has modulus $m_i$ for $i=1,2$.  Then exactly one of the following holds:
\begin{itemize}
    \item[(a)] $\centL \simeq_k k[t_1^{\pm 1},t_2^{\pm 2}]$.
    \item[(b)] $\centL \simeq_k \strange$ for some $\rho\in k^\times$.
\end{itemize}
Moreover (a) holds if and only if $z_1^{m_1}z_2^j \in \bC(\cL)$ for some $j\in \bbZ$
(see Remark \ref{rem:centcalc}).
\end{proposition}

\begin{definition}
\label{def:kind} As in Proposition \ref{prop:kind}, let $\cL =
L(\cA,\Sigma_1,\Sigma_2)$ be a $2$-step iterated loop algebra
based on a finite dimensional central simple algebra $\cA$ over
$k$, where $\Sigma_i$ has modulus $m_i$ for $i=1,2$.  We say that
$\cL$  is of the \emph{first kind} (resp.~\emph{second kind}) if
$\centL$ is isomorphic to $k[t_1^{\pm 1},t_2^{\pm 2}]$
(resp.~$\strange$  for some $\rho\in k^\times$).
\end{definition}

\begin{remark}

(a) It follows from Corollary \ref{cor:centmulti} that
any $2$-step multiloop algebra
based on a finite dimensional central simple
algebra is of  the first kind.

(b) Suppose that $k$ is algebraically closed of characteristic 0
and $\cL = L(\cA,\sigma_1,\sigma_2)$ is a 2-step iterated loop algebra
based on a finite dimensional central simple Lie algebra
$\cA$ over $k$,
where $\sigma_i$ has period $m_i$ for $i=1,2$.
Then $L(\cA,\sigma_1)$ is the derived algebra modulo
its centre of an affine Kac-Moody Lie algebra $\mathfrak g$ \cite[Theorem 8.5]{K2}.  Moreover
one can show that the 2-step loop algebra $\cL$ is of  the first kind in the sense
of Definition \ref{def:kind} if and only if
the automorphism $\sigma_2$ of $L(\cA,\sigma_1)$ is
induced by an automorphism of   the first kind of $\mathfrak g$
(as defined for example  in
\cite[Part III.1]{L}).  Indeed this example is the reason for our choice of terminology.
\end{remark}

We conclude by looking at two examples of
2-step iterated loop algebras.  These examples
illustrate the above proposition (Proposition \ref{prop:kind})
as well as a number of the concepts studied in this article.

\begin{example}
\label{ex:hermitian}
Suppose that $k$ is of characteristic $0$.
In this example we consider a $2$-step iterated loop algebra
$\cL = L(\cA,\sigma_1,\sigma_2)$ based on the Lie algebra
$\cA = \sll$ over $k$, where $\ell \ge 1$ and $\sigma_1$ and $\sigma_2$
have order $m_1=m_2=2$.

Before beginning it will be convenient to define four commuting automorphisms
$\eta_1$, $\eta_2$, $\kappa_1$ and $\kappa_2$ of $\Stwo$
by
\begin{gather*}
\eta_1(z_1^{i_1}z_2^{i_2}) = (-1)^{i_1}z_1^{i_1}z_2^{i_2}, \quad
\eta_2(z_1^{i_1}z_2^{i_2}) = (-1)^{i_2}z_1^{i_1}z_2^{i_2},\\
\kappa_1(z_1^{i_1}z_2^{i_2}) = z_1^{-i_1}z_2^{i_2} \andd
\kappa_2(z_1^{i_1}z_2^{i_2}) = z_1^{i_1}z_2^{-i_2}
\end{gather*}
for $i_1,i_2\in \bbZ$.
Each of these automorphisms restricts to an automorphism of $k[z_1^{\pm 1}]$
which we also denote by $\eta_1$, $\eta_2$, $\kappa_1$ and $\kappa_2$ respectively.

To  construct $\cL$ we first let
$\cL_0 = \cA$.  Next
let $\sigma_1\in \Aut(\cA)$
be defined by $\sigma_1(a) = -Ja^t J$, where
\[J = \begin{bmatrix} 0&\dots&1\\ \vdots&\iddots&\vdots \\ 1 & \dots & 0 \end{bmatrix}.\]
Then $\sigma_1$ has order 2 and we set
\[\cL_1 := L(\cA,\sigma_1,z_1),\]
using the notation of Remark \ref{rem:primitive}.  Thus $\cL_1$ is the algebra
of fixed points in $\cA\otk k[z_1^{\pm 1}]$ of the automorphism $\sigma_1\ot \eta_1$.
(If $k$ is an
algebraically closed field of characteristic 0  and $\ell \ge 2$,
then $\cL_1$ is the derived algebra modulo its centre of the affine Kac-Moody Lie
algebra of type $A_\ell^{(2)}$ \cite[Chapter 8]{K2}.)

Next the automorphisms $1_\cA\ot \kappa_1$ and $\sigma_1\ot
\eta_1$ of $\cA\otk k[z_1^{\pm 1}]$ commute, so  $1_\cA\ot
\kappa_1$ stabilizes $\cL_1$.  We set $\sigma_2 = 1_\cA\ot
\kappa_1 \mid_{\cL_1} \in \Aut_k(\cL_1)$.  Then $\sigma_2$ has
order 2, and we set
\[\cL = \cL_2 := L(\cL_1,\sigma_2,z_2).\]
By construction $\cL$ is a 2-step iterated loop algebra based on $\cA$.

It is clear from the above descriptions of $\cL_1$ and $\cL_2$, that
$\cL$ is the algebra of common fixed points in $\cA\otk \Stwo$ of the automorphisms
$\sigma_1\ot \eta_1$ and $1_\cA \ot \kappa_1\eta_2$.  From this it follows easily
that
\begin{equation}
\label{eq:hermitian}
\cL = \set{x\in \speciallinear_{\ell+1}(K) \suchthat x^* = -x},
\end{equation}
where
\[K = (\Stwo)^{\kappa_1\eta_2}\]
is the algebra of fixed points in $\Stwo$ of the automorphism $\kappa_1\eta_2$, and
\begin{equation}
\label{eq:stardef}
x^* = -J(\eta_1 x)^t J
\end{equation}
for $x\in M_n(K)$. (Here $\eta_1 x$ denotes
the matrix obtained from $x$ by applying $\eta_1$ to the entries of $x$.)
In more geometric language, $\cL$ can
be viewed as the Lie algebra of $K$-linear transformations of the free $K$-module
$K^{\ell +1}$ that are skew relative to the hermitian form $(u,v) \mapsto (\eta_1 u)^t J v$.

Now by Remark \ref{rem:centcalc},
the centroid  of $\cL$
is isomorphic to the algebra
\begin{equation}
\label{eq:Cbarcalc}
\bC(\cL) = \set{ u\in \Stwo \suchthat u\cdot \cL \subset \cL}
\end{equation}
of $\Stwo$.
This together with \eqref{eq:hermitian}
implies that $\bC(\cL)\subset K$.
But by \eqref{eq:stardef},
$(u\cdot x)^* = (\eta_1 u)\cdot x^*$
for $u\in K$ and $x\in \speciallinear_{\ell+1}(K)$. Hence it follows from
\eqref{eq:hermitian} and \eqref{eq:Cbarcalc}
that $\bC(\cL)  = K^{\eta_1}$. So we have
\[\bC(\cL)  = (\Stwo)^{\langle \eta_1 ,\kappa_1\eta_2\rangle}.
\]

Note also that, by Theorem \ref{thm:untwist},  $\Stwo$ is a free
$\bC(\cL)$-module of rank~$4$  and
\[\cL\ot_{\bC(\cL)} \Stwo \simeq  \speciallinear_{\ell+1}(\Stwo).\]
Moreover, by Corollary \ref{cor:centsimplebaseb}, $\cL$ is a prime
pfgc Lie algebra of   type $A_\ell$ (see
Example~\ref{ex:typeLie2}).

Finally, note that $\kappa_1\eta_2(z_1^2z_2^j) = (-1)^jz_1^{-2}z_2^j \ne z_1^2z_2^j$
and so $z_1^2z_2^j \notin \bC(\cL)$ for $j\in \bbZ$.
Thus $\cL$ is of the second kind.  (In fact one can check directly that $\bC(\cL)$ is isomorphic
to $k[u_1,u_2^{\pm 1},w]_\rho$ for $\rho =1$.) So
$\bC(\cL)$ is not isomorphic to
the algebra of Laurent polynomials in any number of variables
(since  $\bC(\cL)$ is not spanned by its units). Hence,
by Corollary \ref{cor:centmulti},
$\cL$ is not
isomorphic to a  multiloop algebra in any number of steps
based on a finite dimensional central simple Lie  algebra.
\end{example}

\begin{example}
\label{ex:quantum}  Suppose that $\ell\ge 1$ and $k$ is a field
which contains a primitive $\ell^\thsup$-root of unity $\zeta = \zeta_\ell$.
In this example we consider a $2$-step multiloop loop algebra
$\cL = M(\cA,\sigma_1,\sigma_2)$ based on the associative
algebra $\cA = M_\ell(k)$ of $\ell\times\ell$-matrices over $k$,
where $\sigma_1$
and $\sigma_2$
have order $m_1=m_2=\ell$.

First let
\[
a_1 = \begin{bmatrix}
1 & 0 & \dots &0\\
0 & \zeta & \dots &0\\
\vdots  & \vdots   & \ddots & \vdots \\
0 & 0 & \dots & \zeta^{\ell - 1}
\end{bmatrix}
\andd
a_2 = \begin{bmatrix}
0 & 1 &  \dots &0\\
\vdots  & \vdots   & \ddots &  \vdots \\
0 & 0 & \dots & 1\\
1 & 0 & \dots & 0
\end{bmatrix}
\]
in $\cA$.  Then $a_2 a_1 = \zeta a_1 a_2$, $a_1^\ell = a_2^\ell = 1$,
and it is well known
that
\[\set{a_1^{i_1} a_2^{i_2} \suchthat 0\le i_1, i_2 \le \ell -1 } \]
is a basis for $\cA$. (See for example \cite[\S 11]{D}.)

Define $\sigma_i\in \Aut_k(\cA)$ by $\sigma_i(x) = a_i x a_i^{-1}$ for $x\in \cA$, $i=1,2$.
Then $\sigma_i(a_i) = a_i$, $\sigma_1(a_2) = \zeta^{-1} a_2$ and $\sigma_2(a_1) = \zeta a_1$.
Hence $\sigma_1$ and $\sigma_2$ are commuting automorphisms of $\cA$ of order $\ell$.
Let
\[\cL = M(\cA,\sigma_1,\sigma_2)\]
be the multiloop algebra of $\sigma_1$, $\sigma_2$
based on $\cA$ (with $m_1 = m_2 = \ell$).  To calculate $\cL$ explicitly note that
$\sigma_1(a_2^{-i_1}a_1^{i_2}) = \zeta^{i_1} a_2^{-i_1}a_1^{i_2}$
and
$\sigma_2(a_2^{-i_1}a_1^{i_2}) = \zeta^{i_2} a_2^{-i_1}a_1^{i_2}$
for $i_1, i_2 \in \bbZ$.  Thus
$\cA_{\modi_1,\modi_2} = k a_2^{-i_1}a_1^{i_2}$ for $i_1, i_2 \in \bbZ$.
Consequently
\[\cL = \spann_k\set{a_2^{-i_1}a_1^{i_2} \ot z_1^{i_1}z_2^{i_2} \suchthat i_1,i_2\in \bbZ}
= \spann_k\set{ x_1^{i_1}x_2^{i_2} \suchthat i_1,i_2\in \bbZ}, \]
where
\[
x_1 = a_2^{-1}\ot z_1 =
\begin{bmatrix}
0 &  \dots & 0 &z_1\\
z_1 & \dots  &  0& 0\\
\vdots & \ddots &\vdots & \vdots \\
0 &  \dots & z_1 & 0
\end{bmatrix}
\text{ and }
x_2 = a_1\ot z_2 =
\begin{bmatrix}
z_2 & 0 & \dots &0\\
0 & \zeta z_2 & \dots &0\\
\vdots  & \vdots   & \ddots & \vdots \\
0 & 0 & \dots & \zeta^{\ell - 1} z_2
\end{bmatrix}
\]
in $\cL$.  Thus $\cL$ is the subalgebra of $M_\ell(\Stwo)$
generated as an algebra by the matrices
$x_1^{\pm 1}$, $x_2^{\pm 1}$ which satisfy the relations
\begin{equation}
\label{eq:quantrel}
x_i x_i^{-1} = x_i^{-1} x_i = 1 \andd x_2 x_1 = \zeta x_1 x_2.
\end{equation}
It follows that $\cL \simeq k_\bq$, where
$k_\bq$ is the algebra presented by the generators
$x_1,x_2$ subject to the relations \eqref{eq:quantrel}.
This  algebra $k_\bq$, which is called the
\textit{quantum torus} determined by the matrix
$\bq = \left[\begin{smallmatrix} 1 & \zeta \\ \zeta^{-1} & 1 \end{smallmatrix}\right]$,
has arisen in a number of different contexts
(see for
example  \cite{M,MP,BGK,G}).

Note that by Corollary \ref{cor:centmulti},
the centroid  (= centre) of $\cL$
is isomorphic to
$\bC(\cL) = k[t_1^{\pm 1},t_2^{\pm 1}]$, where
$t_1 = z_1^\ell$ and $t_2 = z_2^\ell$.
Moreover, by Theorem 6.1,  $\Stwo$ is a free
$\bC(\cL)$-module of rank $\ell^2$ and
$\cL\ot_{\bC(\cL)} \Stwo \simeq  M_\ell(\Stwo)$.
Consequently (see Corollary \ref{cor:Azumaya})
$\cL \simeq k_\bq$ is a prime Azumaya algebra of constant rank~$\ell^2$
that is split by the extension $\Stwo/k[t_1^{\pm 1},t_2^{\pm 1}]$.
\end{example}

\begin{remark}
The fact that the quantum torus $k_{\bq}$ (described in the preceding example)
is an Azumaya algebra was seen
by  a different method some  time ago in \cite[Lemma 4]{M}.
This  information about the algebra $k_\bq$ is important because it tells us that
$k_{\bq}$ defines an element $[k_{\bq}]$ of the  Brauer group
of the ring $k[t_1^{\pm 1},t_2^{\pm 2}]$.
In fact
${}_\ell\hspace{-1pt}\operatorname{Br}(k[t_1^{\pm 1},t_2^{\pm 2}])$
is cyclic of order $\ell$
and the element $[k_{\bq}]$ is a generator of this group \cite[Theorem~6]{M}.
\end{remark}

\begin{remark}  The authors wish to thank John Faulkner for conversations
that led to Example \ref{ex:quantum}.   This example turns out to be
a special case of a more general
construction of quantum tori and their nonassociative
analogs as multiloop algebras. This topic will be
investigated in a article in preparation by the present authors together with
John  Faulkner.
\end{remark}


\begin{thebibliography}{AABGP}

\bibitem[AABGP]{AABGP}
B. Allison, S.~Azam, S.~Berman, Y.~Gao and A.~Pianzola,
\emph{Extended affine Lie algebras and their root systems},
Mem.~Amer.~Math.~Soc. {\bf 126} (603), 1997.

\bibitem[ABGP]{ABGP}
B.~Allison, S.~Berman, Y.~Gao and A.~Pianzola,
\emph{A characterization of affine Kac Moody Lie algebras},
Comm. Math. Phys. {\bf 185} (1997), 671--688.

\bibitem[ABP1]{ABP1} B.~Allison, S.~Berman and A.~Pianzola,
\emph{Covering algebras I: Extended affine Lie algebras},  J.~Algebra
{\bf 250} (2002), 485--516.

\bibitem[ABP2]{ABP2} B.~Allison, S.~Berman and A.~Pianzola,
\emph{Covering algebras II: Isomorphism of loop algebras},
J. Reine Angew. Math. {\bf 571} (2004), 39--71.

\bibitem[ABP3]{ABP3} B.~Allison, S.~Berman and A.~Pianzola,
\emph{Covering algebras III: The nullity 2 case},
in preparation.

\bibitem[BZ]{BZ}
G.~Benkart and  E.~Zelmanov, \emph{Lie algebras
graded by finite root systems and intersection matrix algebras},
Invent. Math. {\bf 126} (1996), 1--45.

\bibitem[BGK]{BGK} S.~Berman, Y.~Gao,  Y.~Krylyuk, \emph{Quantum tori and
the structure of elliptic quasi-simple Lie algebras},
J. Funct. Anal. {\bf 135} (1996), 339--389.

\bibitem[BM]{BM} S.~Berman and R.~V.~Moody, \emph{Lie algebras graded by
finite root systems and the intersection matrix algebras of Slodowy},
Invent. Math {\bf 108} (1992), 323--347.

\bibitem[B:Alg]{B:Alg} N.~Bourbaki, \emph{Algebra. I. Chapters 1--3},
Translated from the
French, Reprint of the 1974 edition, Elements of Mathematics
(Berlin), Springer-Verlag, Berlin, 1989.

\bibitem[B:CA]{B:CA}
N.~Bourbaki, \emph{Commutative algebra. Chapters 1--7}, Translated
from the French, Reprint of the 1989 English translation, Elements
of Mathematics (Berlin), Springer-Verlag, Berlin, 1989.

\bibitem[D]{D} P.K.~Draxl, \emph{Skew fields},
Cambridge Univ. Press, Cambridge, 1983.

\bibitem[EMO]{EMO} T.S.~Erickson,  W.S.~Martindale, 3rd and J.M.~Osborn,
\emph{Prime nonassociative algebras}, Pacific J. Math. {\bf 60}
(1975), 49--63.

\bibitem[G]{G} Y.~Gao, \emph{Representations of
extended affine Lie algebras coordinatized by
certain quantum tori}, Compositio Math. {\bf 123} (2000),  1--25.

\bibitem[H] {H} S.~Helgason, \emph{Differential Geometry, Lie Groups
and symmetric spaces}, Graduate Studies in Mathematics 34, Amer.
Math. Soc., Providence, R.I., 2001.

\bibitem[J1]{J1} N.~Jacobson, \emph{Lie algebras}, Dover, New York, 1979.

\bibitem[J2]{J2} N.~Jacobson, \emph{Structure and representations of
Jordan algebras}, Amer. Math. Soc. Coll. Publ., XXXIX. Amer. Math.
Soc, Providence, R.I., 1968.

\bibitem[K1]{K1}
V.G.~Kac, \emph{Automorphisms of finite order of semi-simple Lie
algebras}, Funct. Anal. Appl. {\bf 3} (1969), 252--254.

\bibitem[K2]{K2} V.G.~Kac,
\emph{Infinite dimensional Lie algebras}, Third edition, Cambridge
University Press, Cambridge, 1990.

\bibitem[Kn]{Kn} M.-A. Knus, \emph{Quadratic and hermitian forms over rings},
Springer-Verlag, Berlin, 1991.

\bibitem[KO]{KO} M.-A. Knus and M.~Ojanguren,
\emph{Th\'eorie de la descente et alg\`ebres d'Azumaya},
Lecture Notes in Mathematics 389, Springer-Verlag, Berlin, 1974.

\bibitem[Ku]{Ku} E.~Kunz, \emph{Introduction to commutative algebra and
algebraic geometry}, Birkhauser, Boston, MA, 1985.

\bibitem[L]{L} F.~Levstein, \emph{A classification of involutive automorphisms
of an affine Kac-Moody Lie algebra}, J. Algebra {\bf 114} (1988), 489--518.

\bibitem[M]{M} A. R. Magid, \emph{Brauer groups of linear algebraic
groups with characters}, Proc. Amer. Math. Soc. {\bf 71} (1978), 164--168.

\bibitem[MP]{MP} J.C.~McConnell and J.J.~ Pettit,
\emph{Crossed products and multiplicative analogues of Weyl algebras},
J. London Math. Soc. {\bf 38} (1988), 47--55.

\bibitem[MZ]{MZ} K.~McCrimmon and E. Zelmanov,
\emph{The structure of strongly prime quadratic Jordan algebras},
Adv. in Math. {\bf 69} (1988),  133--222.

\bibitem[P]{P} A.~Pianzola,
\emph{Affine Kac-Moody Lie algebras as torsors over the punctured line},
Indagationes Mathematicae N.S. {\bf 13} (2002), 249--257.

\bibitem[PS]{PS} S.V.~Polikarpov and I.P.~Shestakov,
\emph{Nonassociative affine algebras}, Algebra and Logic {\bf 29} (1990),  458--466 (1991).

\bibitem[Po]{Po}  U.~Pollmann,  \emph{Realisation der biaffinen Wurzelsysteme
von Saito in Lie--Algebren}, Hamburger Beitr\"age zur Mathematik
aus dem Mathematischen Seminar, Heft {\bf 29}, 1994.

\bibitem[SY]{SY} K.~Saito and D.~Yoshii, Extended affine root systems IV
(Simply laced elliptic Lie
algebras), Publ. RIMS. Kyoto Univ. {\bf 36} (2000), 385--421.

\bibitem[vdL]{vdL} J.~van de Leur, \emph{Twisted Toroidal Lie Algebras},
ArXiv Mathematics e-prints, e-print math/0106119, 2001.

\bibitem[W]{W}  M.~Wakimoto,  \emph{Extended affine Lie algebras and a
certain series of Hermitian representations}, preprint, 1985.

\bibitem[ZSSS]{ZSSS}
K.A.~Zhevlakov,  A.M.~Slin'ko, I.P.~Shestakov, and A.I.~Shirshov,
\emph{Rings that are nearly associative}, Pure and Applied Mathematics, 104,
Academic Press, Inc., New York-London, 1982.

\end{thebibliography}
\end{document}